\documentclass[12pt,a4paper]{article}
\usepackage{amsmath}
\usepackage{amssymb}
\usepackage{amsthm}
\usepackage{amsfonts}
\usepackage{latexsym}
\usepackage{verbatim} 
\usepackage{xcolor}

\usepackage[flushmargin]{footmisc}
\usepackage{color}                    
\usepackage{esvect}
\usepackage{mathrsfs}

\theoremstyle{plain}
\newtheorem{thm}{Theorem}[section]
\newtheorem{cor}{Corollary}[section]
\newtheorem{lem}{Lemma}[section]
\newtheorem{prop}{Proposition}[section]

\theoremstyle{remark}

\theoremstyle{definition}

\newtheorem{exmp}{Example}[section]

\newtheorem{rem}{Remark}[section]

\title{Scaling limits for equivariant Szeg\H{o} kernels near a singular point of the moment map and applications to K\"ahler reduction}
\author{
	Andrea Galasso\footnote{{\bf e-mail}: andrea.galasso@unipegaso.it, andrea.galasso.91@gmail.com, {\bf Address:} Dipartimento di Ingegneria, Universit\`a Pegaso, Sede Centrale Centro Direzionale Isola F2 - Napoli, Italy, {\bf ORCID ID:} 0000-0002-5792-1674 }}
\date{}

\begin{document}
\maketitle

\begin{abstract}  Suppose that a compact and connected Lie group acts holomorphically on a Hodge manifold, and that the action linearises to the associated Hermitian ample line bundle. Assume that zero is a singular value of the associated moment map. The spaces of global holomorphic sections of powers of the line bundle decompose over the finite dimensional irreducible representations of the group. The main goal of this paper is to study the scaling limits for equivariant Szeg\H{o} kernels near a singular point of the moment map. As an application, we discuss the quantization of the singular reduced space: the behavior of invariant sections along the singular strata, the failure of the asymptotic unitarity of the Guillemin--Sternberg map near the singular strata, and the number of invariant sections concentrated near a singular orbit.
\end{abstract}



\newpage
\tableofcontents
\newpage

\section{Introduction}

Let $(M,\,\omega)$ be a $2n$ dimensional Hodge manifold with complex structure $J$ and Riemannian metric $g$. Let $A\,:\,G\times M \rightarrow M$ be a Hamiltonian and holomorphic action of a compact connected Lie group $G$ of dimension $d$. Let $\Phi\,:\,M \rightarrow \mathfrak{g}^{\vee}$ be the moment map, where $\mathfrak{g}:=\mathrm{Lie}(G)$ is the Lie algebra of $G$. 

Let $(L,\,h^L)$ be the quantizing line bundle. For every $k =1,2,\dots$ there is a natural Hermitian structure on each space of holomorphic global sections $H^0(M,L^k)$, and a naturally induced unitary representation of $G$ on $H^0(M,L^k)$. Let $\hat{G}$ be the unitary dual of $G$, we have
\[H^0(M,L^k)=\bigoplus_{\nu\in \hat{G}}H^0(M,L^k)_{\nu}\,.\]
The aim of this paper is to study the asymptotic expansions of the projector onto the isotypes $H^0(M,L^k)_{\nu}$, for a fixed $\nu$, as $k$ goes to infinity. In particular, we focus on the case when $0$ is not a regular value of the moment map. Here, we prefer to work with an equivalent description of these spaces as described in the next paragraphs.

Let $X\subset L^{\vee}$ be the unit circle bundle, with projection $\pi : X \to M$. Then the action of $G$ on $L^{\vee}$ leaves $X$ invariant. Furthermore, $X$ is a
contact and Cauchy-Riemann manifold, with CR structure induced by pulling back $J$ and contact form given by the connection $1$-form $\alpha$. Since $G$ preserves both the Hermitian metric and the holomorphic structure, it
preserves the unique compatible connection, and therefore it acts on $X$ as a group of CR contactomorphisms via $\tilde{A}\,:\,G\times X \rightarrow X$. In particular, $X$ carries a standard $G$-invariant Riemannian metric. 

As is well known, the spaces of smooth sections of $L^{\otimes k}$ may be unitarily and equivariantly identified with the spaces $C^\infty(X)_k$ of smooth
functions on $X$ of the $k$-th isotype for the circle action. Let $H(X)_k\subset C^\infty(X)_k$ be the subspace corresponding to $H^0(M,L^{\otimes k})$ under this isomorphism. The Hardy space $H(X)$ of $X$ is the closure in $L^2(X)$ of
\[
\bigoplus_{k=0}^{\infty} H(X)_k\,.
\]
Thus the representation of $G$ on $H^0(M,L^{\otimes k})$ translates into
\begin{equation}\label{eq:decomp}
	H(X)_k \;=\; \bigoplus_{\nu\in\hat{G}} H(X)_{{k,\nu}}.
\end{equation}
where $\hat{G}$ is the space of unitary irreducible representations of $G$.

Now we can state more precisely our aim. Fix a
pair $(\nu,k)$ in  $\hat{G}\times \mathbb{N}$. 
Let $\Pi_{{k,\nu}}$ be the distributional kernel of the orthogonal projection
\[
\Pi_{{k,\nu}}\,:\, L^2(X)\longrightarrow H(X)_{{k,\nu}}.
\]
It is known that $\Pi_{{k,\nu}}(x,\,x)$ has rapidly decreasing asymptotic expansion as $k$ goes to infinity when $x$ does not lie in the zero locus of the moment map. Here, we shall study the asymptotic properties of the functions $\Pi_{{k,\nu}}(x,\,x)$, with $\Phi(\pi(x))=0$, as $\nu$ is fixed and $k\to +\infty$ when $0$ is not a regular value. This situation is new since, for example in \cite{pscal} or \cite{mz}, $0$ is assumed to be a regular value. Furthermore, it should be mentioned that the techniques adopted in this paper are based on the method introduced in \cite{pscal} and \cite{sz}. The asymptotics of the Szeg\H{o} kernel near a singular point of the moment map seems to be an important case of study, we refer to Section \ref{sec:pre} for examples.

Recall the description of singular points of a moment map in terms of the group action. Infinitesimally, the definition of the moment map says that for any $\xi \in \mathfrak{g}$, and any vector $\mathbf{v} \in T_mM$, we have
\begin{equation}\label{eq:moment-def}
	\langle\mathrm{d}_m\Phi(\mathbf{v}),\, \xi\rangle \;=\; \omega_m(\xi_M(m), \mathbf{v}).
\end{equation}
Here $\langle\cdot,\cdot\rangle$ denotes the pairing between the Lie algebra $\mathfrak{g}$ and its dual $\mathfrak{g}^{\vee}$, and $\xi_M$ is the vector field on $M$ induced by the action of $\xi$.

Equation \eqref{eq:moment-def} implies that the annihilator
$(\mathrm{Im}\, \mathrm{d}\Phi_m)^{\circ}$ of the image of $\mathrm{d}_m\Phi$ in $\mathfrak{g}^{\vee}$
satisfies
\begin{equation}\label{eq:stab-ann}
	(\mathrm{Im}\, \mathrm{d}_m\Phi)^{\circ} \;=\; \{ \xi \in \mathfrak{g} \;\mid\; \xi_M(m) = 0\}\,.
\end{equation}
The right-hand side of \eqref{eq:stab-ann} is precisely $\mathfrak{g}_m$, the Lie algebra of the stabilizer group $G_m$ of $m$. Consequently, $\mathrm{d}_m\Phi$ is surjective if and only if the stabilizer of $m$ is discrete. Equivalently,
\begin{equation*}
	m \text{ is a singular point of }\Phi
	\quad \Longleftrightarrow \quad
	\dim G_m {\,\geq 1}\,.
\end{equation*}
The non-surjectivity of $\mathrm{d}_m\Phi$ is a first obstacle in studying the asymptotic expansion of $\Pi_{{k,\nu}}$ since the moment map enters in the phase of the oscillatory integral and stationary phase formula can not be applied directly.

To express our result, we recall some basic facts about the local geometry of $M$ along $\Phi^{-1}(0)$, see~\cite{ggk}. If $0 \in \mathfrak{g}^{\vee}$ is a singular value of the moment map $\Phi$, 
then $\Phi^{-1}(0)$ is not a submanifold in general.  But, for any $m \in \Phi^{-1}(0)$, we can introduce a local model for a neighborhood of $m$ in $M$. Let $\mathfrak{g}_M(m) \subseteq T_mM$ denote the tangent space to the orbit through $m$, with say $\dim_{\mathbb{R}}(\mathfrak{g}_M(m))=s$. Then, dimension of $G_m$ as a compact submanifold of $G$ is $d-s$, where $d=\dim_{\mathbb{R}}(G)$. Define $N_m \subseteq T_mM$ as the Riemannian orthocomplement of $\mathfrak{g}_M(m)$ in $\mathfrak{g}_M(m)^{\omega}$, see Section \ref{sec:pre}. Then $N_m$ is a complex and symplectic subspace of $T_mM$, of real dimension $2n - 2s$. In Section \ref{sec:pre}, we show the following orthogonal direct sum decompositions:
\begin{equation}\label{eq:orth-decomp}
	T_mM \;=\; \ker\mathrm d_m\Phi \;\oplus\; J_m\left(\mathfrak{g}_M(m)\right)\,, 
	\qquad 
	\ker\mathrm d_m\Phi \;=\; N_m \;\oplus\; \mathfrak{g}_M(m).
\end{equation}
Given the decompositions in \eqref{eq:orth-decomp}, let $m \in M$ and 
$\mathbf{v} \in T_mM$. We may write
\[
\mathbf{v} \;=\; \mathbf{v}_{\mathrm{v}} + \mathbf{v}_{\mathrm{h}} + \mathbf{v}_{\mathrm{t}}\,,
\]
where 
\[
\mathbf{v}_{\mathrm{v}} \in \mathfrak{g}_M(m), 
\quad 
\mathbf{v}_{\mathrm{h}} \in N_m, 
\quad 
\mathbf{v}_{\mathrm{t}} \in J_m\!\big(\mathfrak{g}_M(m)\big).
\]
Here the subscripts $\mathrm{v}, \mathrm{h}, \mathrm{t}$ stand for \emph{vertical}, 
\emph{horizontal}, and \emph{transverse}, respectively.  

By \cite[Corollary B.47]{ggk}, suppose that $M$ is connected, then there exists a closed
subgroup $H \subset G$ such that the set $M_{(H)}$ is open and dense in $M$. To avoid degenerate cases, let us now assume that the action of $G$ is locally free at some point of $M$ (hence locally free generically). In general, for any
subgroup $H\subseteq G$, its fixed-point set
\[
\mathrm{Fix}(H) \;=\; \{ m \in M \;\mid\; h\cdot m = m \text{ for all } h\in H\}
\]
is a union of symplectic submanifolds (possibly of different dimensions). Therefore the critical points of $\Phi$ are a union of symplectic submanifolds of $M$ of the form $\mathrm{Fix}(H)$, where $H$ is a subgroup of $G$, and the singular values of $\Phi$ are the images of these submanifolds. 

The decomposition \eqref{eq:orth-decomp} reflects the following geometric situation. Let $M^m_{(H)}$ be the connected component of the $(H)$ orbit type manifold that contains a given element $m\in M$ such that $\Phi(m)=0$. Then, by \cite[Theorem 8.1.1 (i)]{or}, the set $M^m_{(H)}\cap \Phi^{-1}(0)$ is a submanifold of $M$. The sets $(\Phi^{-1}(0)\cap M_{(H)})/G$, if not empty, are disjoint smooth manifolds, and their union over $H\subset G$ is a partition of $\Phi^{-1}(0)/G$. In this picture, 
$\mathfrak{g}_M(m)$ corresponds to the vertical tangent directions (the fibers of the fibration), while $N_m$ contains the horizontal directions: its $G_m$-fixed part $N_m^{G_m}$ projects isomorphically and unitarily onto $$T_{[m]}\big((\Phi^{-1}(0)\cap M_{(H)})/G\big)\,.$$ Write $N_m=N_m^{G_m}\oplus W_m$, with $W_m$ the $G_m$-invariant complement. The quadratic moment map $\Phi^N$ of the isotropy representation (see \S\ref{sec:pre}) satisfies $$\Phi^N(\mathbf u+\mathbf w)=\Phi^N(\mathbf w)$$ for $\mathbf u\in N_m^{G_m}$, $\mathbf w\in W_m$, and by the local normal form a neighborhood of $[m]$ in $\Phi^{-1}(0)/G$ is homeomorphic to a neighborhood of $(0,0)$ in $$N_m^{G_m}\times\big((\Phi^N|_{W_m})^{-1}(0)/G_m\big),$$ see \cite{sl}. 

Before stating our theorem, another definition is in order. Let $x\in X$, $\pi(x)=m$, and $\Phi(m)=0$. Let $g\in G_{\pi(x)}$, then there exists $h_g\in S^1$ such that $\tilde{A}_g(x)=h_g\cdot x$ since $\pi(\tilde{A}_g(x))={A}_g(\pi(x))=\pi(x)$. Now, we claim that $h_g$ is unique. In fact, say that there exists $h_g'\in S^1$ such that  $\tilde{A}_g(x)=h_g'\cdot x$. Thus, $h_g\cdot x=\tilde{A}_g(x)=h_g'\cdot x$ which implies $h_g'\cdot h_g^{-1}=1$ since the circle action is free on $X$. Then, we shall define
\[A_{\nu,k}(g,x):= 2^{s/2}\,\frac{\dim({\nu})}{V_{\mathrm{eff}}(x)}\cdot \frac{1}{\mathrm{Vol}(G_{\pi(x)})}\,\overline{\chi_\nu(g)\,h_g^k}\qquad (g\in G_{\pi(x)})\,, \]
where $\chi_{\nu}$ is the character corresponding to $\nu$, $\dim({\nu})$ is the dimension of the unitary irreducible representation $\nu$, $V_{\mathrm{eff}}(x)$ is the volume of the orbit in $M$ through $\pi(x)$, and $\mathrm{Vol}(G_{\pi(x)})$ is the volume of $G_{\pi(x)}$ with respect to the Haar measure in $G$.

Here $\mathrm{Vol}(G_{\pi(x)})$ is computed with the Riemannian volume induced by the bi-invariant metric of $G$. Thus $$\mathrm{Vol}(G_{\pi(x)})^{-1}\,\mathrm{dV}_{G_{\pi(x)}}$$ is the normalized Haar measure of $G_{\pi(x)}$, and the final constant does not depend on the metric on $\mathfrak g$. When $s=0$ we set $V_{\mathrm{eff}}(x)=1$. 
	
Since $\Phi(m)=0$, the identity component $G_m^0$ acts trivially on the fibre of $X$ over $m$. Indeed, for $\xi\in\mathfrak g_m$ the generator of the lifted action at $x$ is the horizontal lift of $\xi_M(m)=0$ minus a multiple of $\langle\Phi(m),\xi\rangle\,\partial_\theta=0$. Hence $h_g$ only depends on the class of $g$ in the finite group $G_m/G_m^0$, and $k\mapsto h_g^k$ is periodic. If $G_m$ is connected, $A_{\nu,k}$ does not depend on $k$.


\begin{thm}\label{thm:main} Let $x\in X$ be such that $\Phi\circ \pi(x)=0$, and $\pi(x)=m$. Let us choose a system of $G_m$-adapted Heisenberg local coordinates centered at $x$. For every $\nu\in \hat{G}$ and $\mathbf{v},\,\mathbf{w} \in T_{\pi(x)}M$, the following asymptotic expansion holds as $k$ goes to infinity:
\begin{multline*}
	\Pi_{k,\nu}\Big(x+\tfrac{\mathbf v}{\sqrt k},x+\tfrac{\mathbf w}{\sqrt k}\Big)\sim\Big(\frac k\pi\Big)^{n-\frac s2}e^{Q(\mathbf v_{\mathrm v}+\mathbf v_{\mathrm t},\,\mathbf w_{\mathrm v}+\mathbf w_{\mathrm t})}\\
	\cdot\int_{G_m}A_{\nu,k}(g,x)\,e^{\psi_2(\mathbf v_{g,\mathrm h},\,\mathbf w_{\mathrm h})}\Big(1+\sum_{j\ge1}a_{\nu j}(x,g,\mathbf v,\mathbf w)\,k^{-\frac j2}\Big)\,\mathrm{dV}_{G_m}(g),
\end{multline*}
	where $\mathbf{v}\mapsto \mathbf{v}_g$ is the isotropy action of $g \in G_m$, $a_{\nu j}(x,g,\cdot,\cdot)$'s are polynomials in $\mathbf{v},\,\mathbf{w}$ with coefficients depending on $x$, 
	\[Q\left(\mathbf{v}_{\mathrm{v}}+\mathbf{v}_{\mathrm{t}},\,\mathbf{w}_{\mathrm{v}}+\mathbf{w}_{\mathrm{t}}\right)=-\lVert \mathbf{v}_{\mathrm{t}}\rVert^2-\lVert \mathbf{w}_{\mathrm{t}}\rVert^2+i\,\left[\omega_m\left(\mathbf{v}_{\mathrm{v}},\mathbf{v}_{\mathrm{t}}\right)-\omega_m\left(\mathbf{w}_{\mathrm{v}},\mathbf{w}_{\mathrm{t}}\right)\right]\,, \]
	and
	\[\psi_2(\mathbf{v},\,\mathbf{w}):=\mathbf{v}\cdot \overline{\mathbf{w}}-\frac{1}{2}\left(\lVert \mathbf{v}\rVert^2+\lVert\mathbf{w}\rVert^2\right) \,. \]
	The expansion holds uniformly for $\|\mathbf v\|,\|\mathbf w\|\le C$.
\end{thm} 

The constants refer to the normalization of \cite{sz}, \cite{pscal}, in which $\Pi_k(x,x)=(k/\pi)^n+O(k^{n-1})$ and $$\mathrm{vol}(M)=\frac{\pi^n}{n!}\int_M c_1(L)^n\,;$$ norms and $V_{\mathrm{eff}}$ refer to $g=\omega(\cdot,J\cdot)$. The formula holds at every $m\in\Phi^{-1}(0)$. If $G_m$ is finite, then $s=d$, $$\mathrm{Vol}(G_m)^{-1}\int_{G_m}\cdot\;\mathrm{dV}_{G_m}$$ is the average over $G_m$, and one recovers the scaling limits of \cite{pscal}.

\begin{rem}\label{rem:matching}
Let $y\in(\Phi\circ\pi)^{-1}(0)$ be such that $G_{\pi(y)}$ is finite. Then (see \cite{pscal}, or the case $s=d$ of the Theorem)
\[
\Pi_{k,\nu}(y,y)\sim\Big(\frac k\pi\Big)^{n-d/2}\frac{2^{d/2}\dim(\nu)}{V_{\mathrm{eff}}(y)}\,\frac1{|{G_{\pi(y)}}|}\sum_{g\in {G_{\pi(y)}}}\chi_\nu(g)\,h_g^k \,.
\]
Let $y=x+\mathbf v/\sqrt k$ with $\mathbf v=\mathbf v_{\mathrm h}\in N_m$. Then {$G\cdot\pi(y)$} fibers over $G\cdot m$, so that
\[
V_{\mathrm{eff}}(y)=V_{\mathrm{eff}}(m)\,\mathrm{Vol}(G_m\cdot\mathbf v)\,k^{-(d-s)/2}\big(1+O(k^{-1/2})\big),
\]
where $\mathrm{Vol}(G_m\cdot\mathbf v)$ is the $(d-s)$-dimensional volume of the $G_m$-orbit of $\mathbf v$ in $(N_m,g_m)$. Formally inserting this into the regular formula gives
\begin{equation}\label{eq:formal}
\Big(\frac k\pi\Big)^{n-\frac s2}\frac{2^{s/2}\dim(\nu)}{V_{\mathrm{eff}}(m)}\cdot\frac{(2\pi)^{\frac{d-s}2}}{\mathrm{Vol}(G_m\cdot\mathbf v)}\cdot\frac1{|{G_{\pi(y)}}|}\sum_{g\in {G_{\pi(y)}}}\chi_\nu(g)h_g^k ,
\end{equation}
that is, the same power $k^{n-s/2}$ as in the Theorem. 

Take now $\mathbf v=\mathbf w=\mathbf v_{\mathrm h}\in N_m$ with $\Phi^N(\mathbf v_{\mathrm h})=0$ and $(G_m)_{\mathbf v_{\mathrm h}}$ finite, so that $x+\mathbf v_{\mathrm h}/\sqrt k$ approaches the regular part of $(\Phi\circ\pi)^{-1}(0)$. Since $G_m$ acts unitarily on $N_m$, $$\Re\,\psi_2(g^{-1}\mathbf v_{\mathrm h},\mathbf v_{\mathrm h})=-\frac12\|g^{-1}\mathbf v_{\mathrm h}-\mathbf v_{\mathrm h}\|^2\,.$$ So, as $\|\mathbf v_{\mathrm h}\|\to+\infty$, the integral over $G_m$ concentrates on $(G_m)_{\mathbf v_{\mathrm h}}$. Near $g_0\in(G_m)_{\mathbf v_{\mathrm h}}$ write $g=g_0e^{\zeta}$, $\zeta\in\mathfrak g_m$. In adapted coordinates $$\psi_2(\mathbf a,\mathbf b)=-\frac12\|\mathbf a-\mathbf b\|^2-i\,\omega_m(\mathbf a,\mathbf b)\,,$$ and $\omega_m(\zeta^2\mathbf u,\mathbf u)=-\omega_m(\zeta\mathbf u,\zeta\mathbf u)=0$; hence
\[
\psi_2(g^{-1}\mathbf v_{\mathrm h},\mathbf v_{\mathrm h})=-\tfrac12\|\zeta\cdot\mathbf v_{\mathrm h}\|^2+2i\,\langle\Phi^N(\mathbf v_{\mathrm h}),\zeta\rangle+O\big(\|\zeta\|^3\|\mathbf v_{\mathrm h}\|^2\big),
\]
where $\langle\Phi^N(\mathbf u),\zeta\rangle=\frac12\omega_m(\zeta\cdot\mathbf u,\mathbf u)$. The linear term vanishes on the null cone. Using $$\mathrm{Vol}(G_m\cdot\mathbf v_{\mathrm h})=\mathrm{Vol}(G_m)\,|\det(\zeta\mapsto\zeta\cdot\mathbf v_{\mathrm h})|/|(G_m)_{\mathbf v_{\mathrm h}}|,$$ Laplace's method gives
\begin{multline*}
\frac1{\mathrm{Vol}(G_m)}\int_{G_m}{\overline{\chi_\nu(g)h_g^k}}\,e^{\psi_2(g^{-1}\mathbf v_{\mathrm h},\mathbf v_{\mathrm h})}\,\mathrm{dV}_{G_m}(g)\\
\sim\frac{(2\pi)^{(d-s)/2}}{\mathrm{Vol}(G_m\cdot\mathbf v_{\mathrm h})}\cdot\frac1{|(G_m)_{\mathbf v_{\mathrm h}}|}\sum_{g\in(G_m)_{\mathbf v_{\mathrm h}}}\chi_\nu(g)h_g^k .
\end{multline*}
Since {$G_{\pi(y)}$} is conjugate to $(G_m)_{\mathbf v_{\mathrm h}}$, the leading term of the theorem coincides with \eqref{eq:formal}. If $\Phi^N(\mathbf v_{\mathrm h})\neq0$, the phase $2i\langle\Phi^N(\mathbf v_{\mathrm h}),\zeta\rangle$ produces instead a Gaussian decay in $\|\Phi^N(\mathbf v_{\mathrm h})\|/\|\mathbf v_{\mathrm h}\|$. This is the scaling-limit counterpart of the rapid decay of $\Pi_{k,\nu}$ away from $\Phi^{-1}(0)$. The model case is worked out in Example \ref{ex:bessel}. 
\end{rem}

The leading term of the theorem has a simple interpretation, which is developed in \S \ref{sec:red}. Since $N_m/\!\!/G_m$ is the local model of the reduced space $M_0=\Phi^{-1}(0)/G$ at $[m]$ \cite{sl}, this is a local form of the principle that quantization commutes with reduction, at a singular value. In Section \ref{sec:red} we also recall the stratification of $M_0$ and the corresponding (costratified) quantum space \cite{sj}, \cite{hue}, and we derive some consequences of the theorem, such us the non-unitarity of the Guillemin--Sternberg map \cite{pao}, \cite{mz}, \cite{hk}, \cite{li}. For circle actions, the global counterpart of our theorem is the Riemann--Roch formula of Delarue, Ioos and Ramacher \cite{dir}.


Several proofs of the fact that quantization commutes with reduction are asymptotic. They show that the Guillemin--Sternberg map $\sigma_k$ is an isomorphism only for $k$ sufficiently large; see for instance \cite[Theorem~0.9]{mz} and \cite{hsiaohuang} (where $0$ is assumed to be a regular value), where this follows from the asymptotics of the $G$-invariant Szeg\H{o} (or Bergman) kernel. In the final part of this introduction, we take the opportunity to observe that such a statement automatically holds in every degree if certain algebraic conditions are satisfied, see Proposition \ref{prop:everydegree}. It only uses the multiplicative structure and the normality of the ring of sections, and it does not require $0$ to be a regular value. In our K\"ahler setting the conclusion is of course known by the work of Guillemin--Sternberg \cite{GS:82} and Sjamaar \cite{sj}. The point is that it follows from any asymptotic version of the statement. Similar ideas can be applied in other settings where an asymptotic version of the principle ``quantization commutes with reduction'' is known, see for instance \cite{gh}, provided that the relevant spaces form a graded ring and the reduction map is multiplicative; this is not the case, for instance, for the spectral spaces of a general Toeplitz operator.

Recall that an element of $H(X)_k$ is a CR function $s$ on $X$ such that $s(r_\theta x)=e^{ik\theta}s(x)$, i.e. a holomorphic section of $L^{k}$, so
\[
R:=\bigoplus_{k\ge0}H(X)_k,\qquad R^G:=\bigoplus_{k\ge0}H(X)_{k,0},
\]
(finite sums) are graded $\mathbb C$-algebras. On the reduced side, let $$Z=(\Phi\circ\pi)^{-1}(0),$$ $X_0=Z/G$ and $p:X_0\to M_0$ be the projection. Let $X_{\mathrm{top}}:=p^{-1}(M_{\mathrm{top}})$, where $M_{\mathrm{top}}$ is the open stratum of $M_0$ (see \S\ref{sec:strat}). Define $H(X_0)_k$ to be the space of continuous functions $f$ on $X_0$ such that $f(r_\theta y)=e^{ik\theta}f(y)$ and $f|_{X_{\mathrm{top}}}$ is CR. Equivalently, $f$ is a continuous section of $L_0^k$ over $M_0$ which is holomorphic on $M_{\mathrm{top}}$. Since $M_0$ is a normal complex space \cite{sj}, these are the holomorphic sections of the orbi-bundle $L_0^k$ on the projective variety $M_0$. Set
\[
R_0:=\bigoplus_{k\ge0}H(X_0)_k.
\]
For $s\in H(X)_{k,0}$, the function $s|_Z$ is continuous, $G$-invariant and of type $k$, so it descends to a continuous function $\sigma_k(s)$ on $X_0$, which is CR on $X_{\mathrm{top}}$. Hence $\sigma_k(s)\in H(X_0)_k$, and
\[
\sigma:=\bigoplus_{k\ge0}\sigma_k:\ R^G\longrightarrow R_0,\qquad \sigma(st)=\sigma(s)\,\sigma(t),
\]
is a homomorphism of graded rings, because the restriction is multiplicative. 

\begin{lem}\label{lem:rings}
	$R^G$ and $R_0$ are integral domains, and $R^G$ is integrally closed in its field of fractions.
\end{lem}

\begin{prop}\label{prop:everydegree}
	Suppose that there exists $k_0$ such that $$\sigma_k:H(X)_{k,0}\to H(X_0)_k$$ is an isomorphism for every $k\ge k_0$. Then $\sigma_k$ is an isomorphism for every $k\ge0$.
\end{prop}

\begin{proof}
	If $k=0$, then $H(X)_{0,0}=\mathbb C$, because $M$ is compact and connected, and $H(X_0)_0=\mathbb C$, because $M_0$ is a compact connected normal complex space \cite{sj}. Let $k\ge1$ and choose $N\ge1$ with $Nk\ge k_0$.
	
	First we prove injectivity. Let $s\in H(X)_{k,0}$ with $\sigma_k(s)=0$. Then $s^N\in H(X)_{Nk,0}$ and $\sigma_{Nk}(s^N)=\sigma_k(s)^N=0$, hence $s^N=0$, and therefore $s=0$.
	
	Let us now prove the surjectivity of $\sigma_k$. Let $0\neq f\in H(X_0)_k$. By hypothesis there are $a\in H(X)_{(N+1)k,0}$ and $b\in H(X)_{Nk,0}$ such that
	\[
	\sigma(a)=f^{N+1},\qquad \sigma(b)=f^N .
	\]
	Since $R_0$ is an integral domain (Lemma \ref{lem:rings}), $f^N\neq0$, hence $b\neq0$. Let $\varphi:=a/b\in\operatorname{Frac}(R^G)$. The elements $a^N$ and $b^{N+1}$ lie in $H(X)_{N(N+1)k,0}$ and
	\[
	\sigma(a^N)=f^{N(N+1)}=\sigma(b^{N+1}).
	\]
	Since $N(N+1)k\ge k_0$, injectivity in that degree gives $a^N=b^{N+1}$, i.e. $\varphi^N=b$. Thus $\varphi$ is a root of the monic polynomial $T^N-b\in R^G[T]$. Since $R^G$ is integrally closed in its field of fractions (Lemma \ref{lem:rings}), we get $\varphi\in R^G$. Moreover $\varphi b=a$, where $a$ and $b\neq0$ are homogeneous. Since $R^G$ is a graded domain, comparing the components of highest and lowest degree shows that $\varphi$ is homogeneous, of degree $(N+1)k-Nk=k$. Hence $\varphi\in H(X)_{k,0}$. Finally
	\[
	\sigma(\varphi)\,f^N=\sigma(\varphi)\,\sigma(b)=\sigma(a)=f^{N+1},
	\]
	so that $\big(\sigma(\varphi)-f\big)f^N=0$ in $R_0$. Since $R_0$ is an integral domain and $f^N\neq0$, we conclude that $\sigma_k(\varphi)=f$.
\end{proof}

\section{Preliminaries}
\label{sec:pre}

\subsection{Examples}

Let $K$ be a compact connected Lie group and let $T$ be a maximal torus of $K$. Associated to the couple $(K,T)$, we define the Weyl group $W$ and a system of positive roots $\Delta^+$, we refer to \cite[Section 5]{gusj} for further details. 


Consider a co-adjoint orbit $\mathcal{O}_{\nu}$ associated to a weight $\nu$ in $\mathrm{Lie}(K)^{\vee}$. Then, $\mathcal{O}_{\nu}$ is a symplectic manifold and we denote by $\omega_{KKS}$ the Kostant-Kirillov-Souriau symplectic form. The action of $T$ on $\mathcal{O}_{\nu}$ is Hamiltonian and the associated moment map is denoted by $\Phi\,:\,\mathcal{O}_{\nu}\rightarrow \mathrm{Lie}(T)^{\vee}$. 

Now, let $\nu_1,\,\nu_2,\,\dots,\,\nu_N$ be all the fundamental weights and their conjugates, and  let $\mathrm{Stab}(\nu_i)$ be their stabilizer groups for the co-adjoint
action of {$K$}. The group $T$ is also a maximal torus of $\mathrm{Stab}(\nu_i)$. Denote the Weyl group of the pair $(\mathrm{Stab}(\nu_i), T)$ by $W_i$. 

\begin{thm}[Theorem 5.2.1 in \cite{gusj}] \label{thm:5.2.1}
	The singular points of the moment map $\Phi\,:\,\mathcal{O}_{\nu}\rightarrow \mathrm{Lie}(T)^{\vee}$ are the symplectic submanifolds
	\[\mathrm{Stab}_K(\nu_i)\cdot w\,\nu \]
	where $w\in W$ and $i=1,\dots,N$. Here $\mathrm{Stab}_K(\nu_i)$ means the stabilizer for the action of $K$. The singular values of $\Phi$ are the convex polytopes
	\[ \text{convex hull of }W_i\cdot w\,\nu  \] 
	$w\in W$, and $i=1,\dots, N$.
\end{thm}

The previous theorem should clarify why we are interested in studying singular points of the moment map. In the setting of the introduction we consider $M=\mathcal{O}_{\nu}$ and $G=T$. Even if $0$ does not lie in the image of the moment map nor $0$ is one of its singular values, we shall explain how our result can be applied thanks to the so-called shifting trick, we refer, for example, to \cite[Section 2]{gal}.

{For $G=T$ a torus, the coadjoint orbit $\mathcal O_\nu=\{\nu\}$ is a point, and the shifting trick reduces to twisting the linearization of $L$ by the character $-\nu$ (for $\nu$ integral). This changes neither $X$ nor its Hardy space; it replaces $\Phi$ by $\Phi-\nu$ and identifies $H(X)_{k,k\nu}$ with the $T$-invariant part of level $k$ of the twisted case, so that $\Pi_{k,k\nu}$ equals the invariant projector of the twisted action. The Theorem, with $\nu$ replaced by the trivial representation, then applies at the points of $\Phi^{-1}(\nu)$. For non-abelian $G$ one works instead on $M\times\mathcal O_\nu$, where $\mathcal O_\nu$ is a coadjoint orbit of $G$, see \cite[Section 2]{gal}.}

   
\begin{exmp}\label{ex:bessel}
 		Let $M=\mathbb{CP}^2$, $L=\mathcal O(1)$, and let $S^1$ act by $$e^{i\vartheta}\cdot[z_0:z_1:z_2]=[z_0:e^{i\vartheta}z_1:e^{-i\vartheta}z_2]\,.$$ Then $\Phi$ is a positive multiple of $$(|z_1|^2-|z_2|^2)/|z|^2$$ and $m=[1:0:0]\in\Phi^{-1}(0)$ is a fixed point: $s=0$, $d=1$, $G_m=S^1$, $N_m=T_mM\cong\mathbb C^2$ with weights $(1,-1)$. The action is free on a punctured neighborhood of $m$. For $\nu=0$ the Theorem gives
 		\begin{align*}
 			\Pi_{k,0}\Big(x+\tfrac{\mathbf v}{\sqrt k},x+\tfrac{\mathbf v}{\sqrt k}\Big)&\sim\Big(\frac k\pi\Big)^{2}\frac1{2\pi}\int_0^{2\pi}e^{(e^{-i\vartheta}-1)|v_1|^2+(e^{i\vartheta}-1)|v_2|^2}\,\mathrm d\vartheta\\
 			&=\Big(\frac k\pi\Big)^{2}e^{-\|\mathbf v\|^2}I_0\big(2|v_1||v_2|\big),
 		\end{align*}
 		where $I_0$ is the modified Bessel function. 
 	\end{exmp}

{
\begin{exmp}\label{ex:polygons}
Let $\boldsymbol{r}=(r_1,\dots,r_N)$ be positive integers and $M=(\mathbb{CP}^1)^N$, with $$L=\mathcal O(r_1)\boxtimes\cdots\boxtimes\mathcal O(r_N),$$ so that $n=N$. Let $\omega_{FS}$ be the Fubini--Study form on $\mathbb{CP}^1$, normalized as in \cite{sz}, so that it is the curvature form of $\mathcal O(1)$. Since $\mathcal O(r)=\mathcal O(1)^{\otimes r}$, the curvature form of $L$ is
\[
\omega=\sum_{j=1}^Nr_j\,\mathrm{pr}_j^*\,\omega_{FS},
\]
where $\mathrm{pr}_j:M\to\mathbb{CP}^1$ is the projection onto the $j$-th factor. Set $G=SU(2)$ and let $G$ act diagonally on $M$,  $$g\cdot(u_1,\dots,u_N)=(g\cdot u_1,\dots,g\cdot u_N),$$ where $SU(2)$ acts on each factor $\mathbb{CP}^1=\mathbb P(\mathbb C^2)$ through its standard representation on $\mathbb C^2$; this action lifts to $L$. Up to a positive constant, $$\Phi(u_1,\dots,u_N)=\sum_{j=1}^Nr_ju_j$$ with $u_j\in S^2\subset\mathbb R^3\cong\mathfrak{su}(2)^\vee$. Thus $\Phi^{-1}(0)/SU(2)$ is the space of closed polygons in $\mathbb R^3$ with side lengths $r_j$ \cite{km}, and
\[
H(X)_{k,0}\cong\mathrm{Inv}_{SU(2)}\big(V_{kr_1}\otimes\cdots\otimes V_{kr_N}\big),
\]
where $V_a$ is the irreducible representation of dimension $a+1$. For generic $\boldsymbol{r}$ the polygon space is smooth, and its quantization was studied by Charles \cite{ch}; here we look at the non-generic case.

A closed polygon $u$ which is not contained in a line has $G_u=\{\pm e\}$, where $e$ is the identity of $G$. Now suppose that $u$ is contained in a line, so that $u_j=\varepsilon_j\mathbf n$ for a unit vector $\mathbf n\in S^2$ and signs $\varepsilon_j=\pm1$ with $$\sum_{j=1}^N\varepsilon_jr_j=0.$$ Such polygons are called \emph{degenerate} \cite{km}. Then $G_u\cong U(1)$ is the maximal torus of $SU(2)$ covering the rotations about the axis $\mathbf n$.

Here $s=2$ and $G_u$ is connected, so $h_g\equiv1$. The orbit of $u$ is $G\cdot u=\{(\varepsilon_jp)_j:\ p\in S^2\}\cong S^2$, because $g\cdot u=(\varepsilon_j\,g\mathbf n)_j$ and $g\mathbf n$ ranges over $S^2$. Each component of the embedding $p\mapsto(\varepsilon_1p,\dots,\varepsilon_Np)$ is either the identity or the antipodal map, both of which are isometries of $\mathbb{CP}^1\cong S^2$. Hence the metric induced on $G\cdot u$ by the K\"ahler metric of $\omega$ is $R$ times the Fubini--Study metric, where $R:=\sum_{j=1}^Nr_j$. Since $(\mathbb{CP}^1,\omega_{FS})$ has area $\pi$, we get $V_{\mathrm{eff}}(u)=\mathrm{Vol}(G\cdot u)=\pi R$. Finally, $G_u$ acts on $N_u\cong\mathbb C^{N-2}$ with weights $\pm1$ up to a common factor, which is irrelevant below. Now, write accordingly $\mathbf v_{\mathrm h}=\mathbf v_++\mathbf v_-$, where $\mathbf v_\pm$ come from the sides with $\varepsilon_j=\pm1$, modulo the orbit directions. For $x$ over $u$, the theorem gives
\begin{multline*}
\Pi_{k,0}\Big(x+\tfrac{\mathbf v}{\sqrt k},x+\tfrac{\mathbf v}{\sqrt k}\Big)\\=\Big(\frac k\pi\Big)^{N-1}\frac{2}{\pi R}\,e^{-2\|\mathbf v_{\mathrm t}\|^2-\|\mathbf v_{\mathrm h}\|^2}\,I_0\big(2\|\mathbf v_+\|\,\|\mathbf v_-\|\big)\big(1+O(k^{-1/2})\big).
\end{multline*}
\end{exmp}
}

\subsection{The normal form}


First, we recall the so called tube theorem for proper group actions, for a proof we refer to \cite[Theorem (2.4.1)]{dk}.
 
\begin{thm}[Tube Theorem] \label{thm:tt}
	Let $A\,:\,G\times M \rightarrow M$ be a smooth action of a compact Lie group $G$ on $M$. Then there exists a $G$-invariant open neighborhood $U$ of $m\in M$ such that the $G$-action on $U$ is equivalent to the action of $G$ on $G\times_{G_m} B$, where $B$ is an open $G_m$-invariant neighborhood of $0$ in $T_{m}M / \mathfrak{g}_M(m)$, on which $G_m$ acts linearly via the infinitesimal action modulo $\mathfrak{g}_M(m)$.
\end{thm}

Now, we recall the Marle-Guillemin-Sternberg normal form, see Theorem \ref{thm:normalform} below, for a proof we refer to \cite{or}. From now on, we suppose that $M$ is a quantizable K\"ahler manifold with Hodge form $\omega$. Let $m\in M$ with $\Phi(m) = 0$. Let $R$ be the radical of $\mathfrak{g}_M(m)$, which is the intersection of $\mathfrak{g}_M(m)$ with its symplectic orthogonal. 

The following lemma is well known; we recall the proof for the convenience of the reader. 

\begin{lem} \label{lem:rad} Let $m\in M$ with $\Phi(m)=0$, then $R=\mathfrak{g}_M(m)$.\end{lem}
\begin{proof}
Let $R$ be the radical of $\mathfrak g_M(m)$ with respect to $\omega_m$, i.e.
$$
R=\mathfrak g_M(m)\cap\big(\mathfrak g_M(m)\big)^{\omega},
$$
where
$$
\big(\mathfrak g_M(m)\big)^{\omega}=\{\mathbf{v}\in T_mM:\,\omega_m(\mathbf{v},\mathbf{w})=0 \ \ \forall \mathbf{w}\in\mathfrak g_M(m)\}.
$$

If $\Phi(m)=0$, then for every $\xi,\eta\in\mathfrak g$ the Hamiltonian action identity gives
$$
\omega_m\big(\xi_M(m),\eta_M(m)\big)=\langle\Phi(m),[\xi,\eta]\rangle=0.
$$

Thus any two generators of the tangent space to the orbit at $m$ are $\omega_m$–orthogonal, which means
$$
\mathfrak g_M(m)\subseteq\big(\mathfrak g_M(m)\big)^{\omega},
$$
i.e. $\mathfrak g_M(m)$ is isotropic. Therefore
$$
R=\mathfrak g_M(m)\cap\big(\mathfrak g_M(m)\big)^{\omega}=\mathfrak g_M(m).
$$
\end{proof}

Note that, in Lemma \ref{lem:rad}, if $\mathfrak g$ is abelian, the same identity holds for every $m$, since $[\xi,\eta]=0$ for all $\xi,\eta$.

Let $N_m$ be complement of $\mathfrak{g}_M(m)$ in $\mathfrak{g}_M(m)^{\omega}$,
 \[N_m := \mathfrak g_M(m)^{\perp_g}\cap \mathfrak g_M(m)^{\omega}.
\]
Note that $N_m$ is a symplectic subspace of $(T_mM,\omega_m)$ and thus $N^{\omega}_m$ is symplectic as well. Since $J_m$ is an isomorphism compatible with both $\omega$ and $g$,  we have $ J_m(N_m) = N_m$.  Then, we obtain $T_mM= N^{\omega}_m\oplus N_m$.  A Witt-Artin decomposition at $m$ is a $G_m$-invariant splitting of $T_mM$ as
\begin{equation} 
	T_mM= \mathfrak{g}_M(m)\oplus L_m  \oplus N_m 
\end{equation}
where $L_m$ is the Lagrangian complement to $\mathfrak{g}_M(m)$ in $N^{\omega}_m$, i.e.
\[
N^{\omega}_m = \mathfrak g_M(m) \oplus L_m, \qquad L_m \text{ Lagrangian in } N^{\omega}_m\,.
\]

 Let $\dim \mathfrak{g}_M(m) = s$, $\dim N_m = 2n-2s$, and $\dim T_m M = 2n$. Since $N_m \subset \mathfrak{g}_M(m)^\omega$, we have $\dim N^\omega_m = 2s.$  Now, recall that $(M,\omega,J,g)$ is a Kähler manifold, and thus we can choose the Lagrangian complement $L_m$ as in the following lemma.

\begin{lem} We can choose $L_m= J_m(\mathfrak{g}_M(m))$. \end{lem}
\begin{proof}
Since $\mathfrak g_M(m)$ is isotropic (Lemma \ref{lem:rad}) and $g_m(\mathbf v,J_m\mathbf w)=\omega_m(\mathbf v,\mathbf w)$ up to sign, we have $\mathfrak g_M(m)\perp J_m\mathfrak g_M(m)$. Since $N_m$ is $J_m$-invariant, $N_m^{\omega}=N_m^{\perp}=\mathfrak g_M(m)\oplus J_m\mathfrak g_M(m)$. Finally $J_m\mathfrak g_M(m)$ is isotropic, because $\omega_m(J_m\mathbf v,J_m\mathbf w)=\omega_m(\mathbf v,\mathbf w)$, and it has dimension $s$, so it is a Lagrangian complement of $\mathfrak g_M(m)$ in $N_m^\omega$.
\end{proof}
 

Fix $G$-invariant scalar product $\langle\cdot,\cdot \rangle_{\mathfrak{g}}$ in $\mathfrak{g}$ so that we identify $\mathfrak{g}$ with its dual. The Witt-Artin decomposition can be analyzed further in terms of Tube Theorem \ref{thm:tt} and the properties of $\Phi$. In fact, we have
\begin{itemize}
	\item The subspace $J_m(\mathfrak{g}_M(m))\oplus N_m$ is a slice at $m$ as in the Tube
	Theorem.
	\item $\mathfrak{g}_M(m)^{\omega}=\ker \mathrm{d}_m\Phi$, and therefore $N_m$ can be chosen to be any $G_m$-invariant complement to $\mathfrak{g}_M(m)$ in $\ker \mathrm{d}_m\Phi$.
	\item The map $f:J_m(\mathfrak{g}_M(m))\rightarrow \mathfrak{g}_m^\circ$ defined by
	\[\langle f(\mathbf{v}),\xi\rangle=\omega_m(\xi_M(m),{\mathbf v})\,,\qquad \mathbf{v}\in J_m(\mathfrak{g}_M(m)),\, {\xi \in \mathfrak{g}}\]
	is a $G_m$-equivariant linear isomorphism. Here, $\mathfrak{g}_m^\circ$ is the annihilator of $\mathfrak{g}_m$ in $\mathfrak{g}^{\vee}$.
\end{itemize}

Since by construction, $N$ is a symplectic vector subspace of $T_mM$, call its symplectic form $\omega^N$. Recall also that since the $G$-action on $M$ is Hamiltonian, then $G_m$ acts on $N$ in a Hamiltonian fashion with an associated quadratic homogeneous momentum map $\Phi^N$, see \cite{or}. 


Let $m\in M$ with $\Phi(m)=0$ and consider the direct product
\[
G\times (\mathfrak{g}_m^\circ\times N_m) \,,
\]
where $N_m$ is a $G_m$–invariant symplectic normal space at $m$. Define an action of $G_m$ by
\begin{align*}
	G_m \times (G \times(\mathfrak{g}_m^\circ\times N_m))\to& \,G\times (\mathfrak{g}_m^\circ\times N_m) \\
	(h, (g,\rho, \mathbf{v}))\mapsto &\,(gh^{-1},\, \operatorname{Ad}_{h^{-1}}^*\rho,\,  h\cdot \mathbf{v})\, .
\end{align*}
This action is free, and since $G_m$ is compact, the quotient space is a manifold,
that we denote by
\[
Y= G \times_{G_m}\left(\mathfrak{g}_m^\circ\times N_m\right)\,,
\]
whose elements are the equivalence classes $[g,\rho, \mathbf{v}]$, where two elements in $G\times
(\mathfrak{g}_m^\circ\times N_m)$ are declared equivalent if they belong to the same $G_m$-orbit.

There is also an action of $G$ on $ G \times\left(\mathfrak{g}_m^\circ\times N_m\right)$, given by
\[
g'\cdot (g,\rho, \mathbf{v})=(g'g,\rho, \mathbf{v})\,.
\]
The actions of $G$ and $G_m$ clearly commute, therefore there is an induced $G$-action
on the quotient $Y$ given by
\[
g'\cdot [g,\rho, \mathbf{v}]=[g'g,\rho, \mathbf{v}]\,.
\]

Furthermore, we can endow $Y$ with a differential two-form $\omega^Y$ as follows. First denote by $\pi$ the
projection map
\begin{align*}
	\pi\,:\,G \times\left(\mathfrak{g}_m^\circ\times N_m\right) \rightarrow& \, G \times_{G_m}\left(\mathfrak{g}_m^\circ\times N_m\right) \\
	(g,\rho, \mathbf{v})\mapsto&\, [g,\rho, \mathbf{v}]\, .
\end{align*}
Let us also use the left trivialization of $TG$, given by
\begin{align*}
	G \times \mathfrak{g} \rightarrow& \,TG  \\
	(g, \xi)\mapsto& \,g\,\xi:=\mathrm{d}_eL_g(\xi) \in T_gG\,,
\end{align*}
where $L_g\,:\,G\to G$, $L_g(g')=gg'$ is the left translation. Then, we can
characterize the tangent space to $[g,\rho, \mathbf{v}]\in Y$ as
\[
T_{[g,\rho, \mathbf{v}]}Y=\left\{\gamma_{\xi,\dot{\rho},\dot{\mathbf{v}}}:=\mathrm{d}_{(g,\rho, \mathbf{v})}\pi (g\cdot\xi,\dot{\rho},\dot{\mathbf{v}}) \,:\,\xi\in\mathfrak{g},\,\dot{\rho}\in \mathfrak{g}_m^\circ,\, \dot{\mathbf{v}}\in N_m \right\}\, .
\]
Next, define $\Phi^N:N_m\to \mathfrak{g}_m^{\vee}$ to be the moment map of the isotropy representation of $G_m$ on $(N_m,\,\omega^N)$.
Then we define the two-form $\omega^Y$ by
\begin{align*}
	\omega^Y_{[g,\rho, \mathbf{v}]}\left(\gamma_{\xi_1,\dot{\rho}_1,\dot{\mathbf{v}}_1},\gamma_{\xi_2,\dot{\rho}_2,\dot{\mathbf{v}}_2}\right) 
	=& \;\langle \dot{\rho}_2 +\mathrm{d}_{\mathbf{v}}\Phi^N(\dot{\mathbf{v}}_2),\xi_1\rangle 
	-\langle \dot{\rho}_1 +\mathrm{d}_{\mathbf{v}}\Phi^N(\dot{\mathbf{v}}_1),\xi_2\rangle\\
	&\;+\langle {\rho} +\Phi^N(\mathbf{v}),[\xi_1,\xi_2]\rangle+ \omega_m\!\left(\xi_{1\,M}(m),\xi_{2\,M}(m)\right) \\
	&\; +\omega^N\!\left(\dot{\mathbf{v}}_1,\dot{\mathbf{v}}_2\right)\,.
\end{align*}

This defines a well–defined $G$–invariant symplectic form on $Y$. The $G$–action on $(Y,\omega^Y)$ is Hamiltonian with moment map $\Phi_Y:Y\rightarrow \mathfrak{g}^{\vee}$,
\[
\Phi_Y([g,\rho,\mathbf{v}])=\operatorname{Ad}_{g^{-1}}^{\vee}\big(\rho+\Phi^N(\mathbf{v})\big)\, .
\]

We are now in a position to state the main result of this section, concerning the construction of the Marle–Guillemin–Sternberg normal form.

\begin{thm} \label{thm:normalform}
	Let $m\in M$ with $\Phi(m)=0$, where $0$ is a singular value of the moment map. Let $N$ be a $G_m$-invariant complement to $\mathfrak{g}_M(m)$ in $\ker(\mathrm{d}_m \Phi)$. Then the two-form $\omega^Y$ is symplectic in a neighborhood of $[e,0,0]\in Y$. Furthermore, the $G$-action on $Y$ defined above is Hamiltonian with momentum map
	\[
	\Phi^Y ([g,\rho,\mathbf{v}]) = \operatorname{Ad}^{\vee}_{g^{-1}}\!\left(\rho+\Phi^N(\mathbf{v})\right)\,.
	\]
	
	There exists a local symplectomorphism 
	$$\phi\,:\, Y \rightarrow M
	$$
	which is $G$-equivariant with respect to the action on $Y$ and the original action on $M$, such that 
	\[
	\phi([e,0,0]) = m,\qquad 
	\phi\big(G\times_{G_m}(\{0\}\times\{0\})\big)=G\cdot m\,.
	\]
\end{thm}

Furthermore, $\phi$ can be chosen so that $\mathrm d\phi_{[e,0,0]}$ is the $G_m$-equivariant unitary isomorphism $\iota$ of \eqref{eq:iota} below. Indeed, $\phi$ is obtained from a $G$-equivariant embedding with that differential at $[e,0,0]$, corrected by an equivariant Moser isotopy which is the identity to first order along $G\cdot m$, because the two symplectic forms coincide on $TM|_{G\cdot m}$.

\subsection{Invariant local coordinates}
\label{sec:lc}

In this section we shall briefly recall {the} definition of Heisenberg local coordinates of \cite{sz} at a given $x\in X$ and we make a little equivariant adjustment to the previous construction inspired by \cite{pscal}.

Let $U$ be a sufficiently small open neighborhood of $\pi(x)$ in $M$. Recall that an adapted coordinate system $(z_1,\,\dots,z_n)$ for $M$ centered at $\pi(x)$ is coordinate system defined via the exponential map on a sufficiently small $B_{2n}\left(0,\,\delta\right)\subset {\mathbb{C}^n\cong (T_{\pi(x)}M,J_{\pi(x)})}$ with radius $\delta>0$, $$\exp_{TM}\,:\,B_{2n}\left(0,\,\delta\right)\rightarrow U\subset M$$ that when $\omega_{\pi(x)}$ is expressed in the $z_i$'s  is the standard symplectic structure on $\mathbb{C}^n$. Thus the choice of the $z_i$'s determines a unitary isomorphism between $T_{\pi(x)}M\cong \mathbb{C}^n$.

Furthermore recall that a preferred local frame $e_L\,:\,U\rightarrow L$ for $L$ at $\pi(x)$, in the sense of \cite{sz}, is a holomorphic local section for $L$ in the neighborhood $U$ of $\pi(x)$, satisfying
\begin{equation}\lVert e_L\left(\pi(x)\right) \rVert =1\,,\quad \nabla e_L\left(\pi(x)\right)=0\,,\quad \nabla^2e_L\left(\pi(x)\right)=-\overline{h}_{\pi(x)}\otimes e_L\left(\pi(x)\right)\,,
	\label{eq:el}
\end{equation}
where $\nabla$  is the covariant derivative of the connection, and $\overline{h}=g+i\,\omega$.  The
local holomorphic frame for $L$ uniquely determines a holomorphic dual local frame $e^{\vee}_L$ for $L^{\vee}$, determined by the condition {$\langle e^{\vee}_L,\,e_L\rangle=1$}.

Then, for a sufficiently small $B_{2n}\left(0,\,\delta\right)\subset {\mathbb{C}^n\cong (T_{\pi(x)}M,J_{\pi(x)})}$ with radius $\delta>0$, we define a system of Heisenberg local coordinates for $X$ centered at $x$ by the map
\begin{align} \label{eq:hlc}
	\rho\,:\, (-\pi,\,\pi)\times B_{2n}(0,\,\delta) &\rightarrow X \\
	(\theta,\,\mathbf{w}) &\mapsto e^{i\,\theta}\,\frac{1}{\lVert e_L^{\vee}(\exp_{TM}(\mathbf{w}))\rVert}\,e_L^{\vee}(\exp_{TM}(\mathbf{w}))\,. \notag
\end{align}
If $\mathbf{w} \in T_{\pi(x)}M$, we shall denote by $x+\mathbf{w}$ the point in $X$ with Heisenberg local coordinates $(0,\,\mathbf{w})$. From now on we shall omit $\exp_{TM}$ in expression such as $\exp_{TM}(\mathbf{w})$ and we write simply $\mathbf{w}$.

Suppose $x \in X$, $\Phi\circ \pi(x) =0$. If $G_{\pi(x)}$ is the stabilizer subgroup of $\pi(x)$, for every $g\in G_{\pi(x)}$ there exists a unique $h_g \in S^1$ such that $\tilde{A}_g(x)=h_g \cdot x$, where $\tilde{A}_g\, :\, X \rightarrow X$ is the contactomorphism induced by $g$.  Let $U \subseteq M$ be a $G_m$-invariant open neighborhood of $m$, and suppose that a local holomorphic frame $e^{\vee}_L$ satisfying \eqref{eq:el} has been chosen on $U$. Clearly, for every $g \in G_m$ we
have $g^*(e^{\vee}_L)(m)=\tilde{A}_g \circ e^{\vee}_L(m)$ (recall that $g^*(e^{\vee}_L)=\tilde{A}_g\circ e^{\vee}_L\circ A_{g^{-1}}$). We may then consider the new frame 
\[\widetilde{e}^{\,\vee}_L=\frac{1}{\mathrm{vol}(G_m)}\,\int_{G_m} h^{-1}_g\,g^*(e^{\vee}_L) \,\mathrm{dV}_{G_m}(g)\,,\]
which is well defined, since $G_m$ is compact. Then $\widetilde{e}^{\,\vee}_L(m)=e^{\vee}_L(m)$, and since the metric and the connection are
$G$-invariant, $\widetilde{e}^{\,\vee}_L$ also satisfies \eqref{eq:el}. Moreover, we now have
\begin{equation} \label{eq:ceq}
	g^*(\widetilde{e}^{\,\vee}_L) = h_g \cdot \widetilde{e}^{\,\vee}_L\,, \quad \forall g \in G_m\,.
\end{equation}
\begin{proof}[Proof of equation \eqref{eq:ceq}] Let $l\in G_m$ and consider
	\begin{equation*} 
		l^*(\widetilde{e}^{\,\vee}_L)=\frac{1}{\mathrm{vol}(G_m)}\,\int_{G_m} h^{-1}_g\,{(l\,g)^*}(e^{\vee}_L) \,\mathrm{dV}_{G_m}(g)\,.
	\end{equation*}
Using the fact that $\mathrm{dV}_{G_m}$ is left-invariant, and changing variables {$g\mapsto l^{-1}g$}, we get
	\begin{equation*} 
		l^*(\widetilde{e}^{\,\vee}_L)=\frac{1}{\mathrm{vol}(G_m)}\,\int_{G_m} h^{-1}_{{l^{-1}g}}\,g^*(e^{\vee}_L) \,\mathrm{dV}_{G_m}(g)\,.
	\end{equation*}
Now, since $h^{-1}_{{l^{-1}g}}=h_l\,h^{-1}_{g}$ we get \eqref{eq:ceq}.
\end{proof}

In the following, the underlying preferred local holomorphic frame in the definition of Heisenberg local coordinates will be assumed to satisfy \eqref{eq:ceq}.

Endow $\mathfrak g_m^{\perp}$ with the $\mathrm{Ad}(G_m)$-invariant scalar product {\[\langle\eta,\zeta\rangle_m:=g_m\big(\eta_M(m),\zeta_M(m)\big)\]} and use it to identify $\mathfrak g_m^\circ\cong\mathfrak g_m^\perp$; the map $f$ of \S\ref{sec:pre} then becomes $J_m\zeta_M(m)\mapsto\zeta$. Consider
\begin{equation}\label{eq:iota}
\iota:\ \mathfrak g_m^\perp\oplus\mathfrak g_m^\perp\oplus N_m\longrightarrow T_mM,\qquad (\eta,\zeta,\mathbf u)\longmapsto \eta_M(m)+J_m\zeta_M(m)+\mathbf u .
\end{equation}
Identify $$\mathfrak g_m^\perp\oplus\mathfrak g_m^\perp\cong\mathfrak g_m^\perp\otimes\mathbb C$$ via $(\eta,\zeta)\mapsto\eta+i\zeta$. Then $\iota$ is $\mathbb C$-linear, isometric and $G_m$-equivariant, because $$\mathrm dA_g(\eta_M(m))=(\mathrm{Ad}_g\eta)_M(m)$$ for $g\in G_m$. Moreover $\iota^*\omega_m$ is the standard symplectic form, because $\mathfrak g_M(m)$ is isotropic (Lemma \ref{lem:rad}) and $N_m\perp \mathfrak g_M(m)\oplus J_m\mathfrak g_M(m)$. Hence:
\begin{itemize}
\item[(i)] $\exp_m\circ\,\iota$ is an adapted chart centred at $m$ (with $\omega_m$, $J_m$ and $g_m$ standard), and it is $G_m$-equivariant. Together with the averaged frame satisfying \eqref{eq:ceq}, it defines Heisenberg local coordinates in which,
\[
\tilde A_g\big(x+(\theta,\mathbf u)\big)=x+\big(\theta+\theta_g,\ \mathrm dA_g\mathbf u\big),\qquad g\in G_m,\ \mathbf u\in T_mM,
\]
since isometries fixing $m$ commute with $\exp_m$.
\item[(ii)] The normal form map can be chosen with $\mathrm d\phi_{[e,0,0]}=\iota$, in the coordinates $(\eta,\rho,\mathbf u)\mapsto[\exp_{\mathfrak g}\eta,\rho,\mathbf u]$ of $Y$. Then $\phi$ and $\exp_m\circ\,\iota$ agree up to $O(\|(\eta,\rho,\mathbf u)\|^2)$. 
\end{itemize}
In particular, all the computations below can be carried out in the chart $\exp_m\circ\iota$; the map $\phi$ is only needed to describe $\Phi^{-1}(0)$ and its stratification near $G\cdot m$.



\section{Proof of the main theorem}

Given the direct sum decomposition of $H(X)_k$ into isotypes $H(X)_{k,\,\nu}$, $\nu\in \hat{G}$, we consider the orthogonal projection $\Pi_{k,\,\nu}\,:\,L^2(X)\rightarrow H(X)_{k,\,\nu}$. Thus, for each $x,\,y\in X$
\begin{multline} \label{eq:sz}
	\Pi_{k,\,\nu}(x,\,y)=\dim_{\mathbb{C}}(V_{\nu})\cdot \int_G \chi_{\nu}\left(g^{-1}\right)\,\Pi_{k}\left(\widetilde{A}_{g^{-1}}(x),\,y\right) \mathrm{dV}_G(g) \\
	=\frac{\dim_{{\mathbb{C}}}(V_{\nu})}{2\,\pi}\cdot\int_{-\pi}^{\pi}\int_G e^{-ik\,\theta}\, \chi_{\nu}\left(g^{-1}\right)\,\Pi\left(\widetilde{A}_{g^{-1}}\circ r_{\theta}(x),\,y\right)\, \mathrm{dV}_G(g)\,\mathrm{d}\theta\,.
\end{multline}
where $r_{\theta}$ denotes the circle action on $X$ (here we may take $\mathrm{Vol}(G)=1$; in the final computation we keep track of $\mathrm{Vol}(G)$, to show that the result does not depend on the metric on $\mathfrak g$).

Let $x\in X\setminus (\Phi\circ \pi)^{-1}(0)$, then it is easy to see that $\Pi_{k,\,\nu}$ has a rapidly decreasing asymptotic expansion on a sufficiently small open neighborhood $U$ of $x$. In particular, it follows in a standard manner by iteratively integrating by parts.

Now, the first part of the proof consists in rewriting the integrand in \eqref{eq:sz} more explicitly and introducing appropriate cut-off functions. This procedure is standard; therefore, we will be somewhat sketchy.

From now on, we fix {$x\in (\Phi\circ \pi)^{-1}(0)$, $m=\pi(x)$,} and we consider {the} rescaled displacements $x_{1,k},\, x_{2,k}$ {(with $\mathbf v_1=\mathbf v$, $\mathbf v_2=\mathbf w$ in the notation of the Theorem)}:
\[ x_{j,k}={x}+\frac{\mathbf{v}_j}{\sqrt{k}}\qquad j=1,\,2\,.\]
Then, again, by standard arguments as in \cite{pscal}, the asymptotic behavior of $\Pi_{k,\,\nu}(x_{1,k},\,x_{2,k})$ is unchanged if the integrand in \eqref{eq:sz} is multiplied by $\varrho_k(g)$, where $\varrho_k\,:\,G\rightarrow\mathbb{R}$ is a cut-off function vanishing outside a small shrinking tubular neighborhood $\mathcal{T}_k(G_{\pi(x)})$ of $G_{\pi(x)}$. More precisely, let
\[\mathcal{T}_k(G_{\pi(x)}):=\left\{g\in G\,:\, \mathrm{dist}_G(g,\,h)<2\,k^{\epsilon-1/2},\,h\in G_{\pi(x)}\right\}\,.\]

Then we have
\begin{multline} \label{eq:pichi}
	\Pi_{k,\,\nu}(x_{1,k},\,x_{2,k})\sim\frac{\dim_{{\mathbb{C}}}(V_{\nu})}{2\,\pi}\cdot\int_{-\pi}^{\pi} \int_{\mathcal{T}_k(G_{\pi(x)})} \mathrm{dV}_G(g) \,\mathrm{d}\theta\, \\ \left[e^{-ik\,\theta}\, \chi_{\nu}\left(g^{-1}\right)\,\Pi\left(\widetilde{A}_{g^{-1}}\circ r_{\theta}(x_{1,k}),\,x_{2,k}\right)\,\varrho_k(g)\right]\,,
\end{multline}

	Let $\exp_G : \mathfrak{g} \to G$ be the exponential map, by abuse of notation we shall write $\exp_G (\xi)=e^{\xi}$. Now, let $k$ sufficiently large, so that  $\exp_G$ restricts to a diffeomorphism over $\mathcal{T}_k(G_{\pi(x)})$. Then, we can parameterize it as
	\[g'\in \mathcal{T}_k(G_{\pi(x)}) \text{ as } g'=g\,\exp_G(\xi/\sqrt{k})\,, \quad\xi\in \mathfrak{g}_{\pi(x)}^{\perp}\,. \] 
    Thus, we may express the integration in 
	$\mathrm{dV}_G(g)$ using the exponential chart.  
	
	To this end, fix an orthonormal basis of $\mathfrak{g}_{\pi(x)}^{\perp}$, so as to unitarily 
	identify $\mathfrak{g}_{\pi(x)}^{\perp}$ with $\mathbb{R}^{s}$.  
	Let $\xi$ be the corresponding linear coordinates on $\mathfrak{g}$. We shall denote by $H_{G_{\pi(x)}}(\xi)\, \mathrm{d}\xi$ the local coordinate expression of the Haar measure $\mathrm{dV}_G$ under the exponential chart. 
	
	Then performing the change of variable
	$\xi \mapsto \xi/\sqrt{k}$ we obtain
\begin{multline*} 
	\Pi_{k,\,\nu}(x_{1,k},\,x_{2,k})\sim k^{-s/2}\,\frac{\dim_{{\mathbb{C}}}(V_{\nu})}{2\,\pi}\cdot\int_{-\pi}^{\pi}  \int_{\mathbb{R}^s} \int_{G_{\pi(x)}}  \mathrm{dV}_{{G_{\pi(x)}}}(g)\, \mathrm{d}\xi \,\mathrm{d}\theta\, \\ \left[e^{-ik\,(\theta+\theta_g)}\, \chi_{\nu}\left(e^{-{\xi}/{\sqrt{k}}}g^{-1}\right)\,\Pi\left(\widetilde{A}_{e^{-\xi/\sqrt{k}}g^{-1}}\circ r_{\theta+\theta_g}(x_{1,k}),\,x_{2,k}\right)\,\varrho_k(\xi)\right]\,,
\end{multline*}
where $h_g=e^{i\theta_g}$ for every $g\in G_{\pi(x)}$, and $$\varrho\left(g\,e^{\xi/\sqrt{k}}\right)=\varrho\left(e^{\xi/\sqrt{k}}\right)=:\varrho_k(\xi)$$ for every $k\in \mathbb{N}$ and $g\in G_{\pi(x)}$, since $\varrho$ can be chosen to be $G_{\pi(x)}$-invariant.

Let us introduce an additional cutoff function. Fix {$\delta > 0$} sufficiently small, and let
$\gamma_0 + \gamma_1 = 1$ be a partition of unity on $(-\pi, \pi)$ with 
$\operatorname{supp}(\gamma_0) \subseteq {(-\delta, \delta)}$ and 
$\operatorname{supp}(\gamma_1) \subseteq {(-\pi, -\delta/2) \cup (\delta/2, \pi)}$. 
With this choice, the kernel $\Pi_{k,\,\nu}\!\left(x_{1,k}, x_{2,k}\right)$ splits into two terms, 
\[
\Pi_{k,\,\nu}\!\left(x_{1,k}, x_{2,k}\right) 
= \sum_{j=0}^1 \Pi_{k,\,\nu}\!\left(x_{1,k}, x_{2,k}\right)_{j},
\]
where each $\Pi_{k,\,\nu}\!\left(x_{1,k}, x_{2,k}\right)_{j}$ is given by 
\eqref{eq:pichi} with the integrand multiplied by $\gamma_j$, $j=0,1$. 
In analogy with \cite[Lemma~4.2]{pscal}, one shows that 
$\Pi_{k,\,\nu}\!\left(x_{1,k}, x_{2,k}\right)_{{1}}$ admits a rapidly decreasing 
asymptotic expansion as $k \to +\infty$.

To proceed, let us recall the parametrix for the Szeg\H{o} kernel constructed in~\cite{bs}. 
Up to a smoothing term, which does not contribute to the asymptotic expansion, we can represent $\Pi$ as a complex Fourier integral operator of the form
\begin{equation}\label{eq:parametrix}
	\Pi(x,y) \;=\; \int_{0}^{+\infty} e^{i t \,\psi(x,y)} \, s(x,y,t) \, \mathrm{d}t 
	\qquad (x,y \in X)\,,
\end{equation}
where $\psi$ is the complex phase functions, and the amplitude $s$ is a semiclassical 
symbol admitting an asymptotic expansion
\[
s(x,y,t) \;\sim\; \sum_{j=0}^{+\infty} t^{\,n-j} \, s_j(x,y)\qquad (x,y \in X)\,.
\]
In view of \eqref{eq:pichi}, inserting \eqref{eq:parametrix} into \eqref{eq:pichi}, 
and multiplying the integrand by $\gamma_0$, we obtain
\begin{align} \label{eq:pichibs}
	\Pi_{k,\,\nu}(x_{1,k},\,x_{2,k})\sim&k^{-s/2+1}\,\frac{\dim_{{\mathbb{C}}}(V_{\nu})}{2\,\pi}\cdot\int_0^{+\infty}\int_{-\pi}^{\pi} \int_{\mathbb{R}^s} \int_{G_{\pi(x)}}  \mathrm{dV}_{{G_{\pi(x)}}}(g)\, \mathrm{d}\xi \,\mathrm{d}\theta\,\mathrm{d}t \\ &\left[e^{i\,k\,\Upsilon_k(x,\,\mathbf{v}_1,\,\mathbf{v}_2,\,g,\,\xi,\,\theta,t)}\,\mathcal{A}_k(x,\,\mathbf{v}_1,\,\mathbf{v}_2,\,g,\,\xi,\,\theta,t)\right]\,. \notag
\end{align}
where we have performed the coordinate change $t\mapsto k\,t$, 
\[\Upsilon_k(x,\,\mathbf{v}_1,\,\mathbf{v}_2,\,g,\,\xi,\,\theta,t)=-\theta_g-\theta+t\,\psi\left(\widetilde{A}_{e^{-\xi/\sqrt{k}}g^{-1}}\circ r_{\theta+\theta_g}(x_{1,k}),\,x_{2,k}\right)\]
and
\begin{multline*}
	\mathcal{A}_k(x,\,\mathbf{v}_1,\,\mathbf{v}_2,\,g,\,\xi,\,\theta,t) \\
 =\gamma_0(\theta)\,\chi_{\nu}\left(e^{-{\xi}/{\sqrt{k}}}g^{-1}\right)\,s\left(\widetilde{A}_{e^{-\xi/\sqrt{k}}g^{-1}}\circ r_{\theta+\theta_g}(x_{1,k}),\,x_{2,k},\,k\,t\right)\,\\\quad\cdot\varrho_k\left(\xi \right)\,H_{G_{\pi(x)}}\left(\frac{\xi}{\sqrt{k}}\right)\,.
\end{multline*}

Our next step will be to Taylor expand $\Upsilon_k$ in descending powers of $k^{1/2}$, by relying on (64) and (65) of \cite{sz}. To this end, we shall need the Heisenberg local coordinates of 
\[
\tilde{A}_{e^{-\xi/\sqrt{k}}} \circ r_{e^{i(\theta + \theta_g)}} 
\Big( \tilde{A}_{g^{-1}} \big( x + \tfrac{\mathbf{v}}{\sqrt{k}} \big) \Big)
\quad \text{when } g \in G_{\pi(x)} \,.
\]

Let us consider the isotropy representation
\[
G_{\pi(x)} \longrightarrow \mathrm{GL}(T_{\pi(x)} M), 
\qquad g \longmapsto \mathrm{d}_{\pi(x)} {A}_g .
\]
For every $g\in G_{\pi(x)}$ and $\mathbf{v} \in T_{\pi(x)} M$, let us set
\[
\mathbf{v}_g := \mathrm{d}_{\pi(x)} A_{g^{-1}} (\mathbf{v}) \in T_{\pi(x)} M \,.
\]

In view of our choice of Heisenberg local coordinates, we then have the following lemma, whose proof is similar to \cite[Lemma 4.3]{pscal}.

\begin{lem}
	Let $x \in X$ be such that $\Phi \circ \pi(x)=0$, $\pi(x)=m$, and fix a system of Heisenberg local coordinates centered at $x$. For every $\mathbf{v} \in T_{\pi(x)}M$, $g\in G_m$, and $\xi \in \mathfrak{g}_m^{\perp}$, as $k \to +\infty$, we have
\[
	\tilde{A}_{e^{-\xi/\sqrt{k}}} \circ \left( r_{e^{i\,\theta_g}} 
	\circ \tilde{A}_{g^{-1}} \left( x + \tfrac{\mathbf{v}}{\sqrt{k}} \right) \right)
	=x+\left(\Theta_k(\xi,\,\mathbf{v}),\,\mathbf{V}_k(\xi,\,\mathbf{v})\right)\,.
\]
where
\[\Theta_k(\xi,\,\mathbf{v})=\frac{1}{k}\,\omega_m(\xi_{M}(m),\,\mathbf{v}_g)+O\left(\left\lVert \left(\frac{\xi}{\sqrt{k}},\frac{\mathbf{v}}{\sqrt{k}}\right) \right\rVert^3\right)\,,\]
and
\[\mathbf{V}_k(\xi,\,\mathbf{v})=\frac{1}{\sqrt{k}}\left(\mathbf{v}_g-\xi_{M}(m)\right)+O\left(\left\lVert \left(\frac{\xi}{\sqrt{k}},\frac{\mathbf{v}}{\sqrt{k}}\right) \right\rVert^2\right)\,.\]
\end{lem}

Let us now set $\psi_2(\mathbf{v}_1,\,\mathbf{v}_2)=\mathbf{v}_1\cdot \overline{\mathbf{v}_2}-(\lVert \mathbf{v}_1 \rVert^2+\lVert \mathbf{v}_2 \rVert^2)/2$ where $\mathbf{v}_1,\,\mathbf{v}_2\in \mathbb{C}^n$. By \cite[(63)-(65)]{sz}, we obtain 
\begin{align*}
\Upsilon_k(x,\,\mathbf{v}_1,\,\mathbf{v}_2,\,g,\,\xi,\,\theta,t) &= i\, t\,(1-e^{i\,\theta})-\theta-\theta_g\\  &+\frac{t}{k}\,e^{i\theta}\left[\omega_m(\xi_M(m),\,\mathbf{v}_{1,g})-i\,\psi_2\left(\mathbf{v}_{1,g}-\xi_M(m),\,\mathbf{v}_2\right)\right] \\
&+t\,e^{i\,\theta}\,R_g\left(\frac{\xi_{M}(m)}{\sqrt{k}},\,\frac{\mathbf{v}_1}{\sqrt{k}},\,\frac{\mathbf{v}_2}{\sqrt{k}}\right)
\end{align*}
where $R_g$ is a smooth function vanishing to third order at the origin.

In a similar way as in \cite{pscal}, we are led to define the decomposition:
\[i\,\omega_m(\xi_M(m),\,\mathbf{v}_{1,g})+\psi_2\left(\mathbf{v}_{1,g}-\xi_M(m),\,\mathbf{v}_2\right)= T^g_{\mathrm{h}}+T^g_{\mathrm{t}}+T_{\mathrm{v}}^g+T^g_{\mathrm{tv}}\] 
where
\[T^{g}_{\mathrm{h}}:=\psi_2\left(\mathbf{v}^g_{1,\mathrm{h}},\,\mathbf{v}_{2,\mathrm{h}}\right)\,, \]
\[T^{g}_{\mathrm{t}}:=-\frac{1}{2}\,\lVert \mathbf{v}^g_{1,\mathrm{t}}-\mathbf{v}_{2,\mathrm{t}}\rVert^2\,, \]
\[T^{g}_{\mathrm{v}}:=-\frac{1}{2}\,\lVert \mathbf{v}^g_{1,\mathrm{v}}-\xi_{M}(m)-\mathbf{v}_{2,\mathrm{v}}\rVert^2\,, \]
\begin{align*}T^{g}_{\mathrm{tv}}:=&-i\,\omega_m\left(\mathbf{v}^g_{1,\mathrm{v}}-\xi_{M}(m)-\mathbf{v}_{2,\mathrm{v}},\,\mathbf{v}^g_{1,\mathrm{t}}+\mathbf{v}_{2,\mathrm{t}}\right) \\ &+ i \left[\omega_m(\mathbf v^g_{1,\mathrm v},\mathbf v^g_{1,\mathrm t})-\,\omega_m(\mathbf v_{2,\mathrm v},\mathbf v_{2,\mathrm t})\right]\,.
\end{align*}

Now, we insert in \eqref{eq:pichibs} the asymptotic expansion for the classical symbol	$s$ appearing in the parametrix for $\Pi$, and use a Taylor expansion in $\xi/\sqrt{k}$, $\mathbf{v}_1/\sqrt{k}$ and $\mathbf{v}_2/\sqrt{k}$ in descending powers of $k^{1/2}$, to define
	\begin{equation*}
			\tilde{\mathcal{A}}_k(x,\,\mathbf{v}_1,\,\mathbf{v}_2,\,g,\,\xi,\,\theta,t)  \sim \sum_{l \geq 0} a_{l}^g(x,\,\mathbf{v}_1,\,\mathbf{v}_2,\,\xi,\,\theta,t)\, k^{n-l/2},
	\end{equation*}
	where each coefficient has the form
	\begin{equation*}
		a_{l}^g(x,\,\mathbf{v}_1,\,\mathbf{v}_2,\,\xi,\,\theta,t) = e^{t\, e^{i\theta}\,(T^{g}_{\mathrm{h}}+T^{g}_{\mathrm{t}}+T^{g}_{\mathrm{v}}+T^{g}_{{\mathrm{tv}}})}\, p_{l}^g(x,\,\mathbf{v}_1,\,\mathbf{v}_2,\,\xi,\,\theta,t),
	\end{equation*}
	and each $p_{l}^g$ is a polynomial in $\mathbf{v}_1,\,\mathbf{v}_2$, and $\xi$ with coefficients depending on $x$, $t$, $\theta$ and other parameters. In particular, the leading coefficient is
	\begin{multline} \label{eq:a0g}
		a_{0}^g(x,\,\mathbf{v}_1,\,\mathbf{v}_2,\,\xi,\,\theta,t) = e^{t\, e^{i\theta}\,(T^{g}_{\mathrm{h}}+T^{g}_{\mathrm{t}}+T^{g}_{\mathrm{v}}+T^{g}_{{\mathrm{tv}}})}\, \gamma_0(\theta)\,\chi_{\nu}\left(e^{-{\xi}/{\sqrt{k}}}g^{-1}\right)\,\\
		\cdot\varrho_k\left(\xi\right)s_0(x,x)\,t^n.
	\end{multline}
	Thus, it is easy to see that \eqref{eq:pichibs} can be written as
	\begin{align}  \label{eq:piknu}
			\Pi_{k,\,\nu}(x_{1,k},\,x_{2,k})\sim&k^{-s/2+1}\,\frac{\dim_{{\mathbb{C}}}(V_{\nu})}{2\,\pi}\cdot\int_0^{+\infty}\int_{-\pi}^{\pi} \int_{\mathbb{R}^s} \int_{G_{\pi(x)}}  \mathrm{dV}_{{G_{\pi(x)}}}(g)\, \mathrm{d}\xi \,\mathrm{d}\theta\,\mathrm{d}t \\ &\left[e^{i\,k\,\Psi_g(t,\theta)}\,\tilde{\mathcal{A}}_k(x,\,\mathbf{v}_1,\,\mathbf{v}_2,\,g,\,\xi,\,\theta,t)\right]\,. \notag
	\end{align}
	where
	\begin{equation} \label{eq:psi}
		\Psi_g(t,\theta):=i\,t(1-e^{i\,\theta}) -(\theta+\theta_g)\,.
	\end{equation}
	
	To determine the leading asymptotics of \eqref{eq:piknu}, we proceed as in \cite{pscal}. We first integrate \eqref{eq:a0g} with respect to $\xi$.	The oscillatory integral has phase $\Psi_g$, given in \eqref{eq:psi}, which coincides with that in the scaling asymptotics of non-equivariant Szeg\H{o} kernels \cite{sz}. The phase has non-negative imaginary part and a unique stationary point at $(t,\theta)=(1,0)$. At this point, the Hessian of $\Psi_g$ is non-degenerate. Therefore, by the stationary phase method for complex oscillatory integrals (Theorem 7.7.5 of \cite{h}), and after some routine calculations, we obtain the asymptotic expansion stated in the theorem.

 Set $$\mathbf a:=\mathbf v^g_{1,\mathrm v}-\xi_M(m)-\mathbf v_{2,\mathrm v}$$ and $$\mathbf b:=\mathbf v^g_{1,\mathrm t}+\mathbf v_{2,\mathrm t}.$$ Identify $\mathfrak g_M(m)\cong\mathbb R^s$ and $J_m\mathfrak g_M(m)\cong i\mathbb R^s$ through $\iota$, so that $\omega_m(\mathbf a,\mathbf b)=\mathbf a\cdot\mathbf b$. Since $g\mapsto g\cdot m$ is a fibration $G\to G\cdot m$ with fibre $G_m$, the map $$\mathfrak g_m^\perp\ni\xi\mapsto\xi_M(m)\in\mathfrak g_M(m)$$ has Jacobian
\[
\mathcal J=\frac{V_{\mathrm{eff}}(m)\,\mathrm{Vol}(G_m)}{\mathrm{Vol}(G)} .
\]
At the critical point $(t,\theta)=(1,0)$ one has $$T^g_{\mathrm v}+T^g_{\mathrm{tv}}=-\frac12\|\mathbf a\|^2-i\,\mathbf a\cdot\mathbf b+i\big[\omega_m(\mathbf v^g_{1,\mathrm v},\mathbf v^g_{1,\mathrm t})-\omega_m(\mathbf v_{2,\mathrm v},\mathbf v_{2,\mathrm t})\big],$$ hence
\[
\int_{\mathfrak g_m^\perp}e^{T^g_{\mathrm v}+T^g_{\mathrm{tv}}}\,\mathrm d\xi=\frac{(2\pi)^{s/2}}{\mathcal J}\,e^{-\frac12\|\mathbf b\|^2+i[\omega_m(\mathbf v^g_{1,\mathrm v},\mathbf v^g_{1,\mathrm t})-\omega_m(\mathbf v_{2,\mathrm v},\mathbf v_{2,\mathrm t})]},
\]
and $$-\frac12\|\mathbf b\|^2+T^g_{\mathrm t}=-\|\mathbf v_{1,\mathrm t}\|^2-\|\mathbf v_{2,\mathrm t}\|^2.$$ Now $G_m$ preserves $\omega_m$, $g_m$ and each summand of $$T_mM=\mathfrak g_M(m)\oplus J_m\mathfrak g_M(m)\oplus N_m.$$ So the exponent equals $Q(\mathbf v_{1,\mathrm v}+\mathbf v_{1,\mathrm t},\mathbf v_{2,\mathrm v}+\mathbf v_{2,\mathrm t})$, which does not depend on $g$ and comes out of the integral over $G_m$. For general $(t,\theta)$ the same Gaussian integral, with the quadratic form multiplied by $t\,e^{i\theta}$, gives $(2\pi/(t e^{i\theta}))^{s/2}$ times the corresponding exponential. The stationary phase in $(t,\theta)$ then produces, exactly as in \cite{sz}, the factor $(k/\pi)^n$. Collecting the constants,
\[
\frac{\dim(\nu)}{\mathrm{Vol}(G)}\cdot k^{-s/2}\cdot\frac{(2\pi)^{s/2}}{\mathcal J}\cdot\Big(\frac k\pi\Big)^n=\Big(\frac k\pi\Big)^{n-s/2}\,\frac{2^{s/2}\dim(\nu)}{V_{\mathrm{eff}}(m)\,\mathrm{Vol}(G_m)},
\]
which is the constant appearing in $A_{\nu,k}$.

\bigskip
\section{Final remarks and applications}\label{sec:red}

In this section we interpret the Theorem \ref{thm:main} in terms of the symplectic quotient $M_0:=\Phi^{-1}(0)/G$. First, we note that the scaling limit of $\Pi_{k,0}$ at a singular point $m$ is the Bargmann--Fock quantization of the local model $N_m/\!\!/G_m$ of $M_0$ (see Proposition \ref{prop:fock}). This yields the behavior of invariant sections along the singular strata and the number of invariant sections concentrated near a singular orbit. 

\subsection{The stratified reduced space}\label{sec:strat}
By \cite{sl}, the connected components of the sets $$(\Phi^{-1}(0)\cap M_{(H)})/G\,,$$ 
where $H$ runs over the closed subgroups of $G$, form a stratification of $M_0$; we write $$M_0=\bigsqcup_{l\in\mathcal L}M_l\,.$$ More precisely, $M_0$, with the quotient topology of $\Phi^{-1}(0)\subseteq M$, is a compact Hausdorff space. Each stratum $M_l$ is locally closed in $M_0$ and is a smooth symplectic manifold. Since $M$ is compact, there are finitely many orbit types, so $\mathcal L$ is finite.

If $M_l\cap\overline{M_{l'}}\neq\emptyset$, then $M_l\subseteq\overline{M_{l'}}$. Hence $l\preceq l'$ if and only if $M_l\subseteq\overline{M_{l'}}$ defines a partial order on $\mathcal L$, and $$\overline{M_{l'}}=\bigsqcup_{l\preceq l'}M_l.$$ Furthermore, recall that, if $[m]\in M_l$, a neighborhood of $[m]$ is homeomorphic to a neighborhood of $(0,o)$ in $N_m^{G_m}\times C_m$, where $$C_m:=(\Phi^N|_{W_m})^{-1}(0)/G_m$$ and $o$ is its vertex. Since $\Phi^N$ is homogeneous quadratic, $C_m$ is the cone over the compact stratified space $$\big((\Phi^N|_{W_m})^{-1}(0)\cap S(W_m)\big)/G_m,$$ where $S(W_m)$ is the unit sphere. The stratum through $[m]$ corresponds to $N_m^{G_m}\times\{o\}$; in particular $\dim M_l=\dim N_m^{G_m}$.
There is a unique open stratum $M_{\mathrm{top}}$, it is connected and dense.


In our holomorphic setting, $M_0$ is homeomorphic to the GIT quotient of $M$ by $G^{\mathbb C}$, which is a normal projective variety \cite{sj}, and $M_0$ is a stratified K\"ahler space in the sense of Huebschmann: each stratum carries a K\"ahler structure \cite{hue}.

\begin{exmp}\label{ex:P12}
In Example \ref{ex:bessel}, the invariant ring is generated by $z_0$ and $z_1z_2$, so $M_0$ is the weighted projective line $\mathbb P(1,2)$. There are three strata. The first is the point $[m]$, $m=[1:0:0]$. The second is the point represented by the orbit $$\{[0:1:e^{i\varphi}]\},$$ whose stabilizer is $\{\pm1\}$, it is an orbifold point. The third is the open stratum. Near $[m]$ the local model is $$C_m=\{|v_1|=|v_2|\}/S^1,$$ and $w=v_1v_2$ identifies $C_m$ with $\mathbb C$. As a complex space, $M_0$ is smooth at $[m]$. The reduced K\"ahler structure, however, is singular there. Indeed $$\omega_m=\sum_jr_j\,\mathrm dr_j\wedge\mathrm d\vartheta_j$$ in polar coordinates $v_j=r_je^{i\vartheta_j}$, so the reduced form on $C_m$ is $$\frac12\,\mathrm d|w|\wedge\mathrm d\arg w,$$ which is not smooth at $w=0$, and the reduced metric on $C_m$ is that of a flat cone of total angle $\pi$. 
\end{exmp}

We noe define {the Hardy space of the reduced space.}\label{sec:hardy0}
Let $$Z:=(\Phi\circ\pi)^{-1}(0)\subseteq X\,,$$ $X_0:=Z/G$, and let $p:X_0\to M_0$ be the projection. As observed after the definition of $A_{\nu,k}$, the group $\{h_g:\,g\in G_m\}$ is finite, so the fibre $$p^{-1}([m])\cong S^1/\{h_g\}$$ is a circle. Over each stratum, $X_l:=p^{-1}(M_l)$ is a circle (orbi)bundle with a connection whose curvature is the reduced symplectic form. Hence $X_l$ is a strictly pseudoconvex CR orbi--manifold, in general it is not compact.

Restriction to $Z$ followed by descent defines the Guillemin--Sternberg map
\[
\sigma_k:\ H(X)_{k,0}\longrightarrow C^0(X_0)_k .
\]
By \cite{GS:82} in the regular case and by \cite{sj} in the singular case, $\sigma_k$ is injective and its image is the space of holomorphic sections of the $k$-th power of the descended line bundle on the projective variety $M_0$ (see also \cite{li}). We may therefore define the Hardy space of the reduced space as
\[
H(X_0):=\overline{\bigoplus_{k\ge0}\sigma_k\big(H(X)_{k,0}\big)}\subseteq L^2(X_0,\mu_0),
\]
where $\mu_0$ is the product of the Liouville measure of $M_{\mathrm{top}}$ and $\mathrm d\theta/2\pi$. The lower strata have measure zero, and $\mu_0(X_0)<\infty$. We denote by $\Pi^0$ the corresponding Szeg\H{o} projector and by $\Pi^0_k$ its $k$-th component.

Thus the vector space $H(X_0)_k=\sigma_k(H(X)_{k,0})$ carries two Hermitian structures: the one of $L^2(X_0,\mu_0)$, quantization after reduction, and the one transported from $L^2(X)$ by $\sigma_k$, reduction after quantization.

Hardy spaces of single strata are not canonical, because $M_l$ is in general not compact. For $l\in\mathcal L$ let $Z_l\subseteq Z$ be the preimage of $\overline{M_l}$. Following \cite[Section 3.3]{hrs} and references therein, set
\[
H_{l,k}(X):=\big\{s\in H(X)_{k,0}:\ s|_{Z_l}=0\big\}^{\perp}\subseteq H(X)_{k,0},
\]
where the orthogonal complement is taken in $H(X)_{k,0}$, with respect to the $L^2(X)$ inner product. If $l\preceq l'$ then $Z_l\subseteq Z_{l'}$, hence $H_{l,k}(X)\subseteq H_{l',k}(X)$; moreover $H_{\mathrm{top},k}=H(X)_{k,0}$, since $\sigma_k$ is injective and $M_{\mathrm{top}}$ is dense. The system $(H_{l,k}(X))_{l\in\mathcal L}$, with the orthogonal projections, is a costratified Hilbert space in the sense of \cite[Section 3.3]{hrs}, see also \cite{hue}. In particular $H(X)_{k,0}$ is not the direct sum of spaces attached to the strata, the strata give a filtration indexed by the poset $(\mathcal L,\preceq)$. The invariant Szeg\H{o} kernel describes this filtration.

\begin{lem}\label{lem:span}
For every $l\in\mathcal L$ and $k\ge0$,
\[
H_{l,k}(X)=\mathrm{span}\big\{\Pi_{k,0}(\cdot,x):\ x\in Z_l\big\}.
\]
\end{lem}
\begin{proof}
For $s\in H(X)_{k,0}$ and $x\in X$ we have $s(x)=\langle s,\Pi_{k,0}(\cdot,x)\rangle_{L^2(X)}$. Hence $s$ vanishes on $Z_l$ if and only if it is orthogonal to the span on the right-hand side, and the claim follows since $H(X)_{k,0}$ is finite dimensional.
\end{proof}

\subsection{The local model and applications}\label{sec:fock}
Let $m\in\Phi^{-1}(0)$. Let $\mathfrak F(N_m)$ be the Bargmann--Fock space of $$(N_m,J_m,g_m)\cong\mathbb C^{n-s},$$ which is the space of entire functions $f$ with $$\int|f(\mathbf u)|^2e^{-\|\mathbf u\|^2}\mathrm d\mathbf u<\infty.$$ After multiplication by $e^{-\|\mathbf u\|^2/2}$, the orthogonal projector of $L^2(N_m,\mathrm d\mathbf u)$ onto $\mathfrak F(N_m)$ has kernel $$\pi^{-(n-s)}e^{\psi_2(\mathbf u,\mathbf u')}.$$ The group $G_m$ acts unitarily on $N_m$, hence on $\mathfrak F(N_m)$ by $(g\cdot f)(\mathbf u)=f(g^{-1}\mathbf u)$. Since $g\mapsto h_g$ is a character of $G_m$, the projector of $\mathfrak F(N_m)$ onto the subspace $\mathfrak F(N_m)^{(k)}$ on which $G_m$ acts by the character $h^k$ has kernel
\[
P_m^{(k)}(\mathbf u,\mathbf u'):=\frac{\pi^{-(n-s)}}{\mathrm{Vol}(G_m)}\int_{G_m}\overline{h_g^k}\;e^{\psi_2(g^{-1}\mathbf u,\,\mathbf u')}\,\mathrm{dV}_{G_m}(g).
\]
If $G_m$ is connected, then $h\equiv1$ and $P_m:=P_m^{(k)}$ is the projector onto the invariant Fock space $\mathfrak F(N_m)^{G_m}$. The latter is the completion of the ring $\mathbb C[N_m]^{G_m}$ of invariant polynomials, that is, of the coordinate ring of the affine quotient $N_m/\!\!/G_m^{\mathbb C}$. By the Kempf--Ness theorem \cite{kn}, this affine variety is homeomorphic to $$(\Phi^N)^{-1}(0)/G_m=N_m^{G_m}\times C_m,$$ the local model of $M_0$ at $[m]$ of \S\ref{sec:strat}.

\begin{prop}\label{prop:fock} If $\nu=0$, as $k\to+\infty$
\begin{multline*}
\Pi_{k,0}\Big(x+\tfrac{\mathbf v}{\sqrt k},x+\tfrac{\mathbf w}{\sqrt k}\Big)=\Big(\frac{2k}{\pi}\Big)^{\frac s2}\frac{e^{Q(\mathbf v_{\mathrm v}+\mathbf v_{\mathrm t},\,\mathbf w_{\mathrm v}+\mathbf w_{\mathrm t})}}{V_{\mathrm{eff}}(m)}\;k^{\,n-s}\,P^{(k)}_m(\mathbf v_{\mathrm h},\mathbf w_{\mathrm h})\\+O\big(k^{\,n-\frac s2-\frac12}\big),
\end{multline*}
uniformly for $\|\mathbf v\|,\|\mathbf w\|\le C$. Moreover, with respect to $N_m=N_m^{G_m}\oplus W_m$, $\mathbf u=\mathbf u_0+\mathbf u_1$, we have
\[
P^{(k)}_m(\mathbf u,\mathbf u')=\pi^{-n_l}\,e^{\psi_2(\mathbf u_0,\mathbf u_0')}\cdot P^{(k)}_{W_m}(\mathbf u_1,\mathbf u_1'),
\]
where $n_l:=\dim_{\mathbb C}N_m^{G_m}=\dim_{\mathbb C}M_l$ and $P^{(k)}_{W_m}$ is defined as $P^{(k)}_m$, with $W_m$ in place of $N_m$.
\end{prop}
\begin{proof}The first statement is the Theorem, since\[(k/\pi)^{n-s/2}\,2^{s/2}\,\pi^{n-s}=(2k/\pi)^{s/2}k^{n-s}.\]The second holds because $\psi_2$ is additive with respect to orthogonal decompositions into complex subspaces, and $G_m$ acts trivially on $N_m^{G_m}$.\end{proof}


Note that the factor $$\pi^{-n_l}k^{n_l}e^{\psi_2(\mathbf u_0,\mathbf u_0')}$$ is the universal scaling limit \cite{sz} of the level-$k$ Szeg\H{o} kernel of a K\"ahler manifold of dimension $n_l$, the stratum $M_l$. The factor $P_{W_m}$ is the projector of the cone $C_m$. In this sense, at the scale $k^{-1/2}$ quantization commutes with reduction locally: the scaling limit of the reduced quantum space $H(X)_{k,0}$ at $[m]$ is the quantization $\mathfrak F(N_m)^{G_m}$ of the local model $N_m/\!\!/G_m$ of the reduced space.

One can modify the Guillemin--Sternberg map to a holomorphic map as in \cite{hsiaohuang}, which can be made asymptotically unitary. In the regular case the obstruction is the effective volume $V_{\mathrm{eff}}$ \cite{pao}, \cite{hk}, \cite[Theorem~0.10]{mz}. The metaplectically corrected map is asymptotically unitary \cite{hk}, also without regularity assumptions on $M_0$ \cite{li}. When points of $\Phi^{-1}(0)$ have positive-dimensional stabilizer, the $d$-dimensional orbit volume vanishes on the corresponding strata of $M_0$, and the weighted inner product of \cite{mz} degenerates there. Proposition \ref{prop:fock} describes what happens at scale $k^{-1/2}$ near such a stratum.

\begin{rem}
For a general $\nu$ the same computation gives, in place of $P^{(k)}_m$, the kernel
$$\sum_\tau\frac{\dim(\nu)}{\dim(\tau)}\,[\nu|_{G_m}\otimes h^k:\tau]\;P_{m,\tau}\,,$$ where $\tau$ runs over the irreducible representations of $G_m$, $[\,\cdot:\tau]$ denotes the multiplicity, and $P_{m,\tau}$ is the projector onto the $\tau$-isotype of $\mathfrak F(N_m)$. In particular the leading term at $\mathbf v=\mathbf w=0$ vanishes unless $\nu|_{G_m}\otimes h^k$ contains the trivial representation.
\end{rem}

Let $[m]\in M_l$ and assume for simplicity that $G_m$ is connected. If $\mathbf u\in N_m^{G_m}$, then $\exp_m(\mathbf u/\sqrt k)$ is fixed by $G_m$ and lies within $O(k^{-1})$ of $\Phi^{-1}(0)\cap M_{(G_m)}$, so that $x+\mathbf u/\sqrt k$ is a point over the stratum up to an error which is negligible at the scale $k^{-1/2}$.

\begin{cor}\label{cor:stratum}
For $\mathbf u,\mathbf u'\in N_m^{G_m}$, uniformly on bounded sets,
\[
\Pi_{k,0}\Big(x+\tfrac{\mathbf u}{\sqrt k},x+\tfrac{\mathbf u'}{\sqrt k}\Big)=\frac{(2/\pi)^{\frac s2}}{\pi^{c}\,V_{\mathrm{eff}}(m)}\;k^{\frac s2+c}\cdot\Big(\frac k\pi\Big)^{n_l}e^{\psi_2(\mathbf u,\mathbf u')}\,\big(1+O(k^{-1/2})\big),
\]
where $c:=\dim_{\mathbb C}W_m=n-s-n_l$.
\end{cor}
\begin{proof}Apply Proposition \ref{prop:fock} with $\mathbf v=\mathbf u$ and $\mathbf w=\mathbf u'$; then $Q=0$, and moreover $P_{W_m}(0,0)=\pi^{-c}$.\end{proof}

\begin{exmp}\label{ex:points}
Suppose that $M_l=\{[m]\}$ is a point, $N_m^{G_m}=0$, and let $x\in\pi^{-1}(m)$. By Lemma \ref{lem:span}, $H_{l,k}(X)$ is the line spanned by $\Pi_{k,0}(\cdot,x)$, and by Corollary \ref{cor:stratum}
\[
\big\|\Pi_{k,0}(\cdot,x)\big\|^2=\Pi_{k,0}(x,x)\sim\Big(\frac{2k}\pi\Big)^{\frac s2}\Big(\frac k\pi\Big)^{n-s}\frac1{V_{\mathrm{eff}}(m)}.
\]
Let $$e_{m,k}:=\frac{\Pi_{k,0}(\cdot,x)}{\Pi_{k,0}(x,x)^{1/2}}.$$ If $[m']\neq[m]$ is another point stratum, then
\[
\big|\langle e_{m,k},e_{m',k}\rangle\big|^2=\frac{|\Pi_{k,0}(x,x')|^2}{\Pi_{k,0}(x,x)\,\Pi_{k,0}(x',x')}=O(k^{-\infty}),
\]
because $$\Pi_{k,0}(x,x')=\int_G\Pi_k(\tilde A_{g^{-1}}(x),x')\,\mathrm{dV}_G(g)$$ and $\Pi_k(y,y')=O(k^{-\infty})$ uniformly for $\mathrm{dist}(\pi(y),\pi(y'))\ge\delta>0$. The quantity $|\langle e_{m,k},e_{m',k}\rangle|^2$ is the analogue, in our setting, of the tunneling probabilities between strata considered in \cite{hrs} for a costratified Hilbert space; we have shown that it is $O(k^{-\infty})$ for distinct point strata. In Example \ref{ex:bessel}, $e_{m,k}$ is the normalized section $z_0^k$. In Example \ref{ex:polygons}, the point strata are the degenerate polygons. 
\end{exmp}

\subsection{Counting invariant sections near a singular orbit}\label{sec:count}
Let us measure how many invariant sections are concentrated near a singular orbit. Let $m\in\Phi^{-1}(0)$ with $G_m$ connected, and for $R>0$ let
\[
T_k(R):=G\cdot\exp_m\Big\{\mathbf u_{\mathrm t}+\mathbf u_{\mathrm h}\in J_m\mathfrak g_M(m)\oplus N_m:\ \|\mathbf u_{\mathrm t}\|\le\tfrac R{\sqrt k},\ \|\mathbf u_{\mathrm h}\|\le\tfrac R{\sqrt k}\Big\}
\]
be the corresponding $G$-invariant tube around $G\cdot m$. Let $h_q(m)$ be the number of linearly independent $G_m$-invariant polynomials on $N_m$ which are homogeneous of degree $q$. Let
\[\mathrm P(a,t):=\frac1{\Gamma(a)}\int_0^t\tau^{a-1}e^{-\tau}\,\mathrm d\tau,\]
be the regularized incomplete Gamma function and
\[ \gamma_s(R):=\Big(\frac2\pi\Big)^{\frac s2}\int_{\{\mathbf t\in\mathbb R^s:\,\|\mathbf t\|\le R\}}e^{-2\|\mathbf t\|^2}\mathrm d\mathbf t\]
be the Gaussian on the ball of radius $R$, so that $\gamma_0\equiv1$ and $\gamma_s(R)\to1$ as $R\to+\infty$.

\begin{cor}\label{cor:count}
For every $R>0$, as $k\to+\infty$,
\[
\int_{\pi^{-1}(T_k(R))}\Pi_{k,0}(x,x)\,\mathrm{dV}_X(x)=\mathcal N_m(R)+O(k^{-1/2}),
\]
where
\[
\begin{aligned}\mathcal N_m(R):=&\;\gamma_s(R)\int_{\{\mathbf u\in N_m:\,\|\mathbf u\|\le R\}}P_m(\mathbf u,\mathbf u)\,\mathrm d\mathbf u\\=&\;\gamma_s(R)\sum_{q\ge0}h_q(m)\;\mathrm P\big(q+n-s,\,R^2\big).\end{aligned}
\]
\end{cor}
\begin{proof}
The exponential map $$(g,\mathbf u)\mapsto g\cdot\exp_m(\mathbf u)$$ induces a diffeomorphism of $G\times_{G_m}B_k(R)$ onto $T_k(R)$ for $k\gg0$, where $B_k(R)\subseteq J_m\mathfrak g_M(m)\oplus N_m$ is the set appearing in the definition of $T_k(R)$, a product of two balls of radius $R/\sqrt k$. 

For a $G$-invariant function $f$ this gives $$\int_{T_k(R)}f\,\mathrm{dV}_M=\int_{B_k(R)} f(\exp_m\mathbf u)\,\mathcal D(\mathbf u)\,\mathrm d\mathbf u,$$ where $\mathcal D$ is smooth and $\mathcal D(0)=V_{\mathrm{eff}}(m)$ is the volume of $G\cdot m$. Now, the function $x\mapsto\Pi_{k,0}(x,x)$ is invariant under $G$ and $S^1$. Set $\mathbf u=\mathbf v/\sqrt k$, so that $\mathrm d\mathbf u=k^{-(2n-s)/2}\mathrm d\mathbf v$, and apply Proposition \ref{prop:fock} with $\mathbf v=\mathbf w$, $\mathbf v_{\mathrm v}=0$. The powers of $k$ cancel, $V_{\mathrm{eff}}(m)$ cancels, and we are left with $$(2/\pi)^{s/2}\int_{B_k(R)} e^{-2\|\mathbf v_{\mathrm t}\|^2}P_m(\mathbf v_{\mathrm h},\mathbf v_{\mathrm h})\,\mathrm d\mathbf v_{\mathrm t}\,\mathrm d\mathbf v_{\mathrm h}+O(k^{-1/2}).$$ Choose an orthonormal basis $(p_j)_{j\ge1}$ of $\mathfrak F(N_m)^{G_m}$ consisting of homogeneous polynomials, so that $$P_m(\mathbf u,\mathbf u)=\sum_{j\ge1}|p_j(\mathbf u)|^2e^{-\|\mathbf u\|^2},$$ the series converging locally uniformly. If $p_j$ has degree $q$, then in polar coordinates $$\int_{\|\mathbf u\|\le R}|p_j|^2e^{-\|\mathbf u\|^2}\mathrm d\mathbf u=\mathrm P(q+n-s,R^2).$$ Summing over $j$, by monotone convergence, gives the claim.
\end{proof}

Thus the number of invariant sections living within distance $R/\sqrt k$ of a singular orbit stays bounded as $k\to+\infty$, and its limit is a smoothed counting function of the Hilbert function $q\mapsto h_q(m)$ of the invariant ring of the isotropy representation. 

Let $N:=\dim_{\mathbb C}N_m$ and
	\[
	F_m(\lambda):=\int_{\{\|\mathbf u\|^2\le\lambda\}}P_m(\mathbf u,\mathbf u)\,\mathrm d\mathbf u=\sum_{q\ge0}h_q(m)\,\mathrm P(q+N,\lambda),
	\]
so that $$\mathcal N_m(R)=\gamma_s(R)\,F_m(R^2).$$ The function $F_m$ is a smoothed count of the invariant polynomials. Indeed, $\mathrm P(a,\lambda)$ is the probability that a Poisson variable of mean $\lambda$ is at least $a$. Hence it is close to $1$ for $a\ll\lambda$ and close to $0$ for $a\gg\lambda$, with a transition window of width $O(\sqrt\lambda)$. This reflects the fact that a normalized invariant polynomial of degree $q$ has its Gaussian mass concentrated near the sphere $\|\mathbf u\|^2=q+N$. Therefore,
	\begin{align*}
	F_m(\lambda)&\approx\#\{\text{linearly indep. invariant polynomials of degree }q\le\lambda-N\}\\
	&=\sum_{q\le\lambda-N}h_q(m),
	\end{align*}
	and $\mathcal N_m(R)$ counts, up to the factor $\gamma_s(R)$ coming from the directions tangent to the orbit, the invariant polynomials whose mass lies in the ball of radius $R$.

\begin{exmp}\label{ex:count}
In Example \ref{ex:bessel}, $\mathbb C[N_m]^{G_m}=\mathbb C[v_1v_2]$, so $h_q=1$ for $q$ even and $h_q=0$ for $q$ odd. Since $\mathrm P(a,t)$ is the probability that a Poisson variable of mean $t$ is at least $a$, we get
\[
\mathcal N_m(R)=\sum_{b\ge0}\mathrm P(2b+2,R^2)=\frac{R^2}2-\frac{1-e^{-2R^2}}4 .
\]
\end{exmp}

Let $\Sigma\subseteq\Phi^{-1}(0)$ be a compact connected component of an orbit type stratum, with $G_m$ connected for $m\in\Sigma$, and let $T_k(R;\Sigma)$ be the tube of radius $R/\sqrt k$ around $\Sigma$. The estimates in the proof of Theorem \ref{thm:main} are uniform for $m\in\Sigma$, because $\Sigma$ is compact, the orbit type is constant along $\Sigma$, and the data entering the proof (normal form, Heisenberg coordinates, isotropy representation) depend smoothly on $m$. Arguing as in Corollary \ref{cor:count}, and by Proposition \ref{prop:fock}, one obtains
\[
\int_{\pi^{-1}(T_k(R;\Sigma))}\Pi_{k,0}(x,x)\,\mathrm{dV}_X(x)\sim\Big(\frac k\pi\Big)^{n_l}\int_{\Sigma/G}\mathcal N_{W_m}(R)\,\mathrm{dV}_{\Sigma/G}([m]),
\]
where $n_l=\dim_{\mathbb C}\Sigma/G$ and $\mathcal N_{W_m}(R)$ is defined as $\mathcal N_m(R)$, with $W_m$ in place of $N_m$. Hence a singular stratum $M_l$ influences $\dim H(X)_{k,0}$ at order $k^{\dim_{\mathbb C}M_l}$.

\bigskip
\textbf{Tool and computational resource disclosure.}
The author used an AI assistant (Anthropic's Claude, November 2026) in the final stage of the document in exploratory discussions concerning some applications. The author takes full responsibility for the correctness and adequacy of the arguments and results, and for the accuracy of the attribution of prior work.

\textbf{Acknowledgments.} The author is a member of GNSAGA (Gruppo Nazionale per le Strutture Algebriche, Geometriche e le loro Applicazioni) of INdAM (Istituto Nazionale di Alta Matematica ``Francesco Severi'') and gratefully acknowledges its support.

\end{document}